\pdfoutput=1
\documentclass[letterpaper, 10 pt, conference]{ieeeconf}  

\IEEEoverridecommandlockouts                              

\usepackage{graphicx}
\usepackage{amsmath} 
\usepackage{amssymb}  

\usepackage{amsfonts}
\usepackage{mathrsfs}
\usepackage{dsfont}
\newtheorem{thm}{\bf Theorem}[section]
\newtheorem{cor}[thm]{\bf Corollary}
\newtheorem{lem}[thm]{\bf Lemma}
\newtheorem{rem}[thm]{\bf Remark}
\newtheorem{ass}[thm]{\bf Assumption}
\newtheorem{prop}[thm]{\bf Proposition}
\newtheorem{prob}[thm]{\bf Problem}
\newtheorem{deff}[thm]{\bf Definition}

\newcommand{\R}{\mathbb{R}}

\newcommand{\ee}{\mathbb{E}}
\newcommand{\LL}{\mathcal{L}}
\newcommand{\cj}{\wedge}
\newcommand{\eps}{\varepsilon}
\newcommand{\ball}{\overline{\mathcal{B}}}
\newcommand{\oball}{\mathcal{B}}
\newcommand{\soln}{\Phi_\delta}
\newcommand{\solnp}{\Phi_{\delta'}}

\newcommand{\kl}{\mathcal{KL}}

\newcommand{\ff}{\mathcal{F}}
\newcommand{\pp}{\mathbb{P}}
\newcommand{\eee}{\mathcal{E}}

\newcommand{\ppp}{\mathcal{P}}

\newcommand{\br}{\mathscr{B}}

\newcommand{\supx}{\sup\limits_{X\in\soln(x,W)}}
\newcommand{\infx}{\inf\limits_{X\in\soln(x,W)}}
\newcommand{\init}{\mathcal{X}_0}
\newcommand{\rr}{\mathcal{R}}
\newcommand{\dd}{\mathcal{D}}
\newcommand{\sys}{\mathcal{S}_\delta}
\newcommand{\sysp}{\mathcal{S}_{\delta'}}
\newcommand{\doa}{\mathcal{G}_{\delta}}
\newcommand{\uu}{\mathscr{U}}

\newcommand{\aba}{\oball_R(A)}
\usepackage{color}

\newcommand{\BlackBox}{\rule{1.5ex}{1.5ex}}    
\newcommand{\pfbox}{\hfill\BlackBox\\[2mm]}  
\newenvironment{pf}{\par\noindent{\bf Proof\ }}{\hfill\BlackBox\\[2mm]}

\usepackage{hyperref}
 \hypersetup{
    colorlinks=true,
    linkcolor=blue,
    filecolor=magenta,      
    urlcolor=cyan,
    pdftitle={Overleaf Example},
    pdfpagemode=FullScreen,
    }
    \usepackage[misc]{ifsym}

\title{\LARGE \bf
Sufficient Conditions for Robust Probabilistic Reach-Avoid-Stay Specifications using Stochastic Lyapunov-Barrier Functions 
}

\author{Yiming Meng and Jun Liu
\thanks{This work was supported in part by the NSERC of Canada and the CRC and ERA programs. Y. Meng and J. Liu are with  the Department of Applied Mathematics, University of Waterloo, Waterloo, Ontario, Canada. 
        {\tt\small yiming.meng@uwaterloo.ca,  j.liu@uwaterloo.ca 
}}%
}

\begin{document}

\maketitle
\thispagestyle{empty}
\pagestyle{empty}

\begin{abstract}
Stability and safety are crucial in safety-critical control of dynamical systems. The reach-avoid-stay objectives for deterministic dynamical systems can be effectively handled by formal methods as well as Lyapunov methods with soundness and approximate completeness guarantees. However, for continuous-time stochastic dynamical systems, probabilistic reach-avoid-stay problems are viewed as challenging tasks. Motivated by the recent surge of applications in characterizing safety-critical properties using Lyapunov-barrier functions, we aim to provide a stochastic version for  probabilistic reach-avoid-stay problems in consideration of robustness. To this end, we first establish a  connection between probabilistic stability with safety constraints and  reach-avoid-stay specifications. We then prove that stochastic Lyapunov-barrier functions provide sufficient conditions for the target objectives. We apply Lyapunov-barrier conditions in control synthesis for reach-avoid-stay specifications, and show its effectiveness in a case study. 
\end{abstract}

\section{Introduction}
The reach-avoid-stay property is one of the building
blocks for specifying more complex temporal logic
objectives. Control synthesis of such specifications has received substantial interests in areas such as robotic motion
planning \cite{fribourg2013control, faulwasser2009model,nilsson2017augmented}.  

In the deterministic context, verification and control synthesis problems are achievable via abstraction-based formal methods \cite{belta2017formal}. 
Considering uncertain transition systems with bounded measurable signals, robust abstractions with soundness and approximate completeness provide guarantees for a given specification \cite{liu2017robust,li2020robustly}. Despite  algorithmic improvements on reducing computational complexities \cite{girard2009approximately,liu2017robust,belta2017formal}, 
it still remains a fundamental challenge to overcome the curse of dimensionality in abstraction-based approaches for verification and control synthesis.

On the other hand, Lyapunov-like functions  are able to connect stability and safety attributes with reach-avoid-stay properties in terms of 
characterizing the approximated domain of attraction \cite{liu2020smooth,  romdlony2016stabilization}. Thanks to the fundamental converse theorems of Lyapunov and barrier functions \cite{teel2000smooth,lin1996smooth,liu2020converse}, the theoretical work in  \cite{liu2020smooth}  justifies that a smooth Lyapunov-barrier function is sufficient and necessary (in a slightly weaker sense) for reach-avoid-stay objectives. In terms of  programming, the framework \cite{ames2016control} also shows effectiveness in the control synthesis of safety and stabilization objectives without spatial discretization. The recent work \cite{meng2021control} took advantages of the above and achieved control synthesis for reach-avoid-stay specifications in application to a system that undergoes a Hopf-bifurcation. For such systems with tunable parameters, the abstraction-based algorithms underperform Lyapunov-barrier approaches  due to the difficulties of adjusting the speed of the dynamical flows.

As we usually concern how probability laws distribute the corresponding weak solutions in the stochastic context, 
a proper specification is to specify a probability of sample paths satisfying certain state-space  behaviors. 
As for verification and control synthesis of probabilistic stability-safety type problems, it appears more challenging.  
Authors in \cite{lahijanian2015formal,cauchi2019efficiency,majumdar2020symbolic,dutreix2020verification,dutreix2020specification} applied abstract models, such as interval-valued Markov chain  and bounded-parameter
Markov decision process,  on discrete-time continuous-state stochastic systems to compute an inclusion of the real satisfying probability and synthesize controllers for probabilistic specifications (including probabilistic reachability on an infinite horizon). Works in \cite{reissig2018symbolic,vinod2017scalable,kariotoglou2013approximate} characterized value functions for reachability/reach-avoid problems in discrete-time continuous-state stochastic systems and applied dynamic programming for synthesizing optimal controllers. Authors in \cite{esfahani2016stochastic} developed a weak dynamic programming principle for the value functions of probabilistic reach-avoid specifications in continuous-time continuous-state stochastic systems, which provides compatibility for non-almost-sure probabilistic requirements. 

Since small perturbations should necessarily be taken into account due to reasons such as  modelling uncertainties and measurement errors of the state, robust analysis provides guarantees in a worst-case scenario. 
Despite the current theme of regarding `inaccuracy' from the  computation of probability measures as the `uncertainty'  \cite{lahijanian2015formal,cauchi2019efficiency,majumdar2020symbolic,dutreix2020verification}, to make a closer analogy of the deterministic case, we consider uncertainties as a result of perturbed stochastic systems which create an inclusion of solutions. 
A robust satisfaction of a probabilistic specification in a perturbed stochastic system is then interpreted as follows: the solution process measured in the correspondingly worst but accurate probability law still satisfies the probabilistic specification. The work \cite{subbaraman2015robust} demonstrated the robust Lyapunov-stability for discrete-time stochastic systems with perturbations. The authors in \cite{poveda2015flexible} considered the same type of systems and developed a robust algorithm to guarantees practical convergence to a Nash equilibrium in non-cooperative
games. Continuous-time stochastic differential inclusions are also well studied \cite{kisielewicz2007stochastic,kisielewicz2013stochastic,malinowski2013interrelation}. 

Motivated by the deterministic robust abstractions \cite{liu2017robust,li2020robustly} and the recent comparisons with robust Lyapunov-type characterizations of reach-avoid-stay specifications \cite{meng2021control}, to better understand how these two connect in the stochastic context, this paper formulates stochastic Lyapunov-barrier functions to deal with sufficient conditions for robust probabilistic reach-avoid-stay specifications. 

The rest of this paper is organized as follows. In Section
\ref{sec: pre}, we present the preliminaries for the systems, concepts of solutions, as well as other important definitions. In Section \ref{sec: conn}, we
show the connections between robust probabilistic reach-avoid-stay and stability with safety guarantees. In Section \ref{sec: lb},
we provide sufficient conditions for robust probabilistic reach-avoid-stay satisfactions. In section \ref{sec: app},
a case study on a stochastic Moore-Greitzer model is conducted to 
demonstrate how controllers can be generated based on a control version of stochastic Lyapunov-barrier certificates. The paper is concluded in Section \ref{sec: conclusion}.

\textbf{Notation:} 
We denote the Euclidean space by $\R^n$ for $n>1$.  We denote $\R$ the set of real numbers, and $\R_{\geq 0}$ the set of nonnegative real numbers.  Given $a,b\in\R$, we define $a\wedge b:=\min(a,b)$. 
Let $C_b(\cdot)$ be the space of all bounded continuous functions/functionals $f: (\cdot)\rightarrow\R$. A  continuous and strictly increasing function $\alpha:\,
\R_{\ge 0}\rightarrow\R_{\ge 0}$ is said to belong to class $\mathcal{K}$ if $\alpha(0)=0$. 

The open ball of radius $r$ centered at $x$ is denoted by $\oball_r(x):=\{y\in\R^n: |y-x|< r\}$, where $|\cdot|$ is the Euclidean norm. We also use $\oball_r:=\oball_r(0)$ to represent  open balls centered at $0$. Given two sets $A,B\subseteq\R^n$, the set difference of $B$ and $A$ is defined by $B\setminus A=\{x\in B: x\not\in A\}$.
For a given set $A\subseteq\R^n$, we denote by $A^c$ the complement of the set $A$ (i.e.,  $\R^n\setminus A$); denote by $\bar{A}$ 
(resp. $\partial A$) the closure (resp. boundary) of $A$.  For a closed set $A\subseteq \R^n$ and $x\in\R^n$, we denote the distance from $x$ to $A$ by $|x|_A=\inf_{y\in A}|x-y|$ and $r$-neighborhood of $A$ by $\oball_r(A)=\bigcup_{x\in A}\oball_r(x)$. 


For any 
stochastic processes $\{X_t\}_{t\geq 0}$ we use the shorthand notation $X:=\{X_t\}_{t\geq 0}$. For any stopped process $\{X_{t\cj\tau}\}_{t\geq 0}$, where $\tau$ is a stopping time, we use the shorthand notation $X^\tau$. We denote the Borel $\sigma$-algebra of a set by $\mathscr{B}(\cdot)$ and the space of all probability measures on $\br(\cdot)$ by $\mathscr{M}(\cdot)$. 

\section{Preliminaries}\label{sec: pre}
\subsection{\bf System dynamics}
Consider the following perturbed stochastic differential equation (SDE):
\begin{equation}\label{E: sys}
    dX_t=f(X_t)dt+\xi(t)dt+g(X_t)dW_t,\;\;X_0=x,
\end{equation}
where $\;\xi:\R_{\geq 0}\rightarrow\ball_\delta $ is any measurable point mass signal; $f:\R^n\rightarrow\R^n$ is a nonlinear vector field; $g:\R^n\rightarrow \R^{n\times m}$ is a smooth mapping; $W$ is an $m$-dimensional Wiener process. For future references, we denote systems driven by SDE \eqref{E: sys} by $\sys$, of which the $\delta$ represents the $\delta$-perturbations.
\begin{ass}\label{ass: usual}
We make the standing assumptions on the regularity of the system $\sys$ for the rest of this paper:
\begin{itemize}
    \item[(i)] The mappings $f,g$ satisfy local Lipschitz continuity.
    \item [(ii)] The eigenvalues $\lambda_i[(gg^T)(x)]$ of the matrix $gg^T(x)$ for $i=1,2,\cdots,n$ satisfy
    $$\sup\limits_{x\in\R^n}\min_{i=1,2,\cdots,n}\lambda_i[(gg^T)(x)]>0.$$
    \item[(iii)]There exists a trivial solution $x_e$ for system $\mathcal{S}_0$ such that $f(x_e)=g(x_e)=0$.
\end{itemize}
\end{ass}
\begin{deff}[Characteristic operator]\label{def: generator}
For each $d\in\ball_\delta$,
we denote by $\LL_d$  the  characteristic operator of 
$\mathcal{S}_\delta$ as 
$$\LL_d h(x)=\nabla h(x)\cdot(f(x)+d)+\frac{1}{2}\operatorname{Tr}\left[(gg^T)(x)\cdot h_{xx}(x)\right],  $$
where
$h\in C^2(\R^n)$, $h_{xx}=(h_{x_ix_j})_{n\times n}$, and $\operatorname{Tr}[\;\cdot\;]$ denotes the trace.
\end{deff}

Since we only care about the probabilistic properties in the state space, we consider mostly the weak solutions of the perturbed SDEs. 
\begin{deff}
The system $\sys$ admits a weak solution if there exists a (most likely unknown) filtered probability space $(\Omega^\dagger,\mathscr{F}^\dagger,\{\mathscr{F}^\dagger_t\},  \pp^\dagger)$, where a Wiener process $W$ is defined and a pair $(X,W)$ is adapted, such that $X$ solves the SDE \eqref{E: sys} for any $\xi:\R_{\geq 0}\rightarrow \ball_\delta$.

We denote by $\soln(x,W)$ the set of all weak solutions with  $X_0=x$ a.s. for a given $x\in\R^n$. Likewise, for a given set $K\subseteq\R^n$, let $\soln(K,W)$ denote the set of all weak solutions with any initial distribution  on $(K,\mathscr{B}(K))$.
\end{deff}

\subsection{\bf Canonical space}
We have a Wiener process $W$ defined on some probability space $(\Omega^\dagger,\mathscr{F}^\dagger, \pp^\dagger)$ for each weak solution. We transfer information to the canonical space, which gives us the convenience to study the law of the solution processes as well as the probabilistic behavior in the state space. 
Define $\Omega:=C([0,\infty);\R^n) $ with coordinate process $\mathfrak{X}_t(\omega):=\omega(t)$ for all $t\geq 0$ and all $\omega\in\Omega$. Define $\ff_t:=\sigma\{\mathfrak{X}_s,\;0\leq s\leq t\}$  for each $t\geq 0$, then  the smallest $\sigma$-algebra containing the sets in every $\ff_t$, i.e.  $\ff:=\bigvee_{t\geq 0}\ff_t$, turns out to be same as $\mathscr{B}( \Omega) $. 
For each $X\in\soln(\R^n,W)$, the induced measure (law)  $\ppp^X\in\mathscr{M}( \Omega)$ on $\ff$ is such that $\ppp^X(A)=\pp^\dagger\circ X^{-1}(A)$ for every $A\in\mathscr{B}( \Omega)$. We also denote $\eee^X$ by the associated expectation operator w.r.t. $\ppp^X$.

To emphasize on the uncertainty of laws of a system $\sys$, we prefer to work on the probability spaces $(\Omega,\ff,  \ppp^X)$ for each weak solution $X$ rather than the original $(\Omega^\dagger,\mathscr{F}^\dagger, \pp^\dagger)$. 

\begin{deff}\textbf{(Weak convergence of measures and processes):}
Given any separable metric space $(\mathcal{S},\rho)$,
a sequence of $\{\ppp^n\}$ of $\mathscr{M}(\mathcal{S})$ is said to weakly converge to $\ppp\in\mathscr{M}(\mathcal{S}) $, denoted by $\ppp^n\rightharpoonup\ppp$, if for all  $f\in C_b(\mathcal{S})$ we have
$\lim_{n\rightarrow\infty}\int_{\mathcal{S}}f\;d\ppp^n=\int_{\mathcal{S}}f\;d\ppp. $ 
A sequence $\{X^n\}$ of continuous processes $X^n$ with law $\ppp^n$ is said to weakly converge (on $[0,T]$) to a continuous process $X$ with law $\ppp^X$, denoted by $X^n\rightharpoonup X$,  if for all $f\in C_b(C([0,T];\R^n))$ we have 
$\lim_{n\rightarrow\infty}\eee^{n}[f(X^{n})]=\eee^X[f(X)].$
\end{deff}

\subsection{\bf Strong Markov properties}
Since for each measurable signal $\xi$, the corresponding martingale problem is well posed under Assumption \ref{ass: usual}, by Markovian selection theorems \cite[Theorem 5.19, Chap 4]{ethier2009markov}, the unique solution $\ppp$ to the martingale problem also makes the weak solution $(X,W)$ Markovian.

\subsection{\bf Other definitions}
We first provide definitions for probabilistic set stability given a closed set $A\subseteq\R^n$.

\begin{deff}[Uniform stability in probability]\label{def: u.s.}
The set $A$ is said to be uniformly stable in probability 
(Pr-U.S.) for $\sys$  if for each $\eps\in(0,1)$ there exists $\varphi_\eps\in\mathcal{K}$ such that
\begin{equation}\label{E: prus}
    \inf\limits_{X\in\soln(x,W)}\ppp^X[|X_{t}|_A\leq \varphi_\eps(|x|_A)\;\;\forall t\geq 0]\geq 1-\eps,
\end{equation}
where $x$ is the initial condition. 
\end{deff}

\begin{rem}\label{rem: epsdelta}
Equation
\eqref{E: prus} is equivalent to the following: for any $\eps\in(0,1)$ and $r>0$, there exists an $\eta=\eta(\eps,r)\in(0,r)$ such that 
\begin{equation}
    \inf\limits_{X\in\soln(x,W)}\ppp^X[|X_t|_A\leq r\;\;\forall t\geq s(\omega)]\geq 1-\eps,
\end{equation}
whenever $|X_{s(\omega)}|_A\leq \eta $ for some random time $s(\omega)$. We can simply pick $\eta=\varphi^{-1}_\eps$.
\end{rem}

\begin{deff}[Uniform attractivity in probability]
The set $A$ is said to be uniformly attractive in probability (Pr-U.A.) for $\sys$ if there exists some $\eta>0$ such that, for each $\eps\in(0,1)$, $r>0$, there exists some $T>0$ such that whenever $|x|_A<\eta$,
\begin{equation}
    \inf\limits_{X\in\soln(x,W)}\ppp^X[|X_{t}|_A<r,\;\forall t\geq T]\geq 1-\eps.
\end{equation}


\end{deff}
\begin{deff}\textbf{(Uniformly asymptotic stability in probability):}
The set $A$ is said to be uniformly
asymptotically stable in probability (Pr-U.A.S.)  for $\sys$ if it is Pr-U.S. and Pr-U.A.  for $\sys$.
\end{deff}

Next we introduce several definitions pertinent to probabilistic stability with safety guarantees. To this end, we consider a closed unsafe set $U\subseteq\R^n$. 
\begin{deff}[Work place]\label{def:workspace}
Since the solutions are not generally non-explosive without stability assumptions, a bounded workplace $\rr:=\oball_{\tilde{R}}(x_e)$ with sufficiently large $\tilde{R}>0$ is added as an extra constraint. We name $\dd=\dd(\rr,U):=\rr\cap U^c$.
\end{deff}

\begin{deff}[Explosion and safety]\label{def:explosion}
For any solution $X\in\soln(\R^n,W)$, we define the corresponding explosion time 
$\sigma^*=\sigma^*(\rr):=\inf\{t\geq 0: X_t\in \rr^c\}$ and safety time $\sigma=\sigma(\dd):=\inf\{t\geq 0: X_t\in \dd^c\}$.
\end{deff}
\begin{rem}\label{rem: nonexplos}
Safety is usually the priority in practice. Given safety requirement w.r.t. $\mathcal{D}$ (resp. $\mathcal{R}$), to study conditional probabilistic properties of some process $X$, it is equivalent 
to just working with
the law of $X^{\sigma}$ (resp. $X^{\sigma^*}$). Note that for systems with trivial Pr-U.S. sets, the indicator $\mathds{1}_{\{\sigma^*=\infty\}}\rightarrow 1$ as $R\rightarrow \infty$ and does not render `too much harm' to replace the law of $ X^{\sigma^*}$ by $\ppp^X$. 
\end{rem}

The following theorem verifies a notion of weak compactness of stopped weak solutions of SDE \eqref{E: sys}. 

\begin{prop}\label{thm: weak_compact}
Under the Assumption \ref{ass: usual}, given any compact set $K$,  the set of all stopped process $X^{\sigma^*}$ is nonempty and sequentially weakly compact (w.r.t. the weak convergence) on every filtered probability space $(\Omega,\ff,\{\ff_t\}_{t\in[0,T]})$, where $X\in\bigcup_{x\in K}\soln(x,W)$ (resp. $X\in\soln(K,W)$).  That is, given any sequence of weak solutions $\{X^n\}_{n=1}^\infty$ in the above sense, 
there is a subsequence $\{X^{n_k}\}$, a  process $X\in\bigcup_{x\in K}\soln(x,W)$ (resp. $X\in\soln(K,W)$) such that $(X^{n_k})^{\sigma^*}\rightharpoonup X^{\sigma^*}$.
\end{prop}
\begin{rem}
The conclusion follows immediately by \cite[Theorem 1]{kisielewicz2007stochastic} and \cite[Corollary 1.1, Chap 3]{kisielewicz2013stochastic}.
The proof falls in standard procedures. We can first show that the truncated laws $\{\ppp^{n,{\sigma^*}}\}$ of the stopped processes $\{(X^n)^{\sigma^*}\}$ form a tight family of measures on $(\Omega,\ff, \{\ff_t\}_{t\in [0,T]})$. Then the relatively weak compactness follows since $(X^{n_k})^{\sigma^*}\rightharpoonup X^{\sigma^*}$ if and only if  $\ppp^{n_k,{\sigma^*}}\rightharpoonup \ppp^{\sigma^*}$. The weak closedness comes from compactness of the reachable sets of the stopped processes.
\end{rem}

Now we introduce two closely-related  specifications pertaining to stability and safety issues. 
\begin{deff}[Probabilistic stability with safety]
Given a closed set $U\subseteq\R^n$, let $\dd$ and $\sigma$ be defined as in Definitions \ref{def:workspace} and \ref{def:explosion}, respectively. 
Given $\init,A\subseteq\dd$ and $p\in [0,1]$,  $\sys$ is said to satisfy a probabilistic stability under safety specification w.r.t. $(\init,A,U)$ with probability at least  $p$, denoted by  $(\init,A, U, p)$, if
\begin{enumerate}
    \item[(1)]$A$ is closed and Pr-U.A.S. for $\sys$;
    \item[(2)] For all $X\in\bigcup_{x\in\init}\soln(x,W)$, 
    $$\ppp^X\left[\sigma=\infty\;\text{and}\;\lim\limits_{t\rightarrow\infty}|X_t|_A=0\right]\geq p. $$
\end{enumerate} 
\end{deff}

\begin{deff}
Given $\init,\Gamma\subseteq \dd$. On $(\Omega,\ff)$, for each $X\in\bigcup_{x\in\init}\soln(x,W)$, we define the events
\begin{enumerate}
    \item[(i)] $RS(\init,\Gamma, \dd):=\{\omega: \gamma<\infty\;\text{and}\; X_{t\cj \sigma}\in\Gamma\;\forall t\geq \gamma\}$, where $\gamma:=\inf\{t\geq 0: X_t\in\Gamma\}$ of $X$;
    \item[(ii)] $RAS(\init,\Gamma,\dd):= RS(\init,\Gamma,\dd)\cap \{\sigma=\infty\}$.
\end{enumerate}
\end{deff}

\begin{deff}(\textbf{Probabilistic reach-and-stay and reach-avoid-stay specification}): 
Given a closed set $U\subseteq\R^n$, let $\dd$ and $\sigma$ be defined as in Definitions \ref{def:workspace} and \ref{def:explosion}, respectively. 
Given $\init,\Gamma\subseteq \dd$ 
and $p\in [0,1]$,  $\sys$ is said to satisfy a reach-avoid-stay specification w.r.t. $(\init,\Gamma, U)$ with probability at least  $p$, denoted by  $(\init,\Gamma, U, p)$, if for every $X\in\bigcup_{x\in\init}\soln(x,W)$, we have 
$\ppp^X[RAS(\init,\Gamma,\dd)]\geq p.$
\end{deff}

\section{A Connection to Probabilistic Stability with Safety Guarantee}\label{sec: conn}
\subsection{\bf Probabilistic stability with safety  implies probabilistic reach-avoid-stay}
We first show that if a closed set $A$ is Pr-U.S. for $\sys$, then any weak solutions starting at $x$ from a compact subset of the $p$-domain of attraction is uniformly attracted to $A$ with probability at least $p$.
\begin{prop}\label{prop: uniform}
Suppose that a closed set $A\subset \dd$ is Pr-U.S. for $\sys$. Let $K$ be a compact set and $p\in(0,1)$. Then the following two statements are equivalent:
\begin{enumerate}
    \item[(1)] For any solution $X\in \bigcup_{x\in K}\soln(x,W)$,
    $$\ppp^X\left[\lim\limits_{t\rightarrow\infty}|X_{t\cj\sigma}|_A=0\right]\geq p.$$
    \item[(2)] For every $r>0$, there exists $T=T(r,\eps)$ such that for any $X\in\bigcup_{x\in K}\soln(x,W)$,
    $$\ppp^X\left[|X_{t\cj \sigma}|_A<r,\;\forall t\geq T\right]\geq p. $$
\end{enumerate}
\end{prop}
\begin{pf}
Clearly (2) implies (1). We only show the converse. Suppose that (2) is not true. 
Then there exists some $r>0$ such that for all $n>0$ there exists $x_n\in K$, $X^n\in\soln(x_n,W)$ with law $\ppp^n$ such that
\begin{equation}\label{E: oppo}
    \ppp^n[|X_{t\cj\sigma}^n|_A\leq r,\;\forall t\geq n]<p.
\end{equation}
 Now let $\tau^n=\inf\{t\geq 0: X_{t\cj\sigma}^n\in \oball_\eta(A)\}$, where $\eta$ is to be chosen later  for each $n$. Rearranging \eqref{E: oppo} we have,
\begin{equation}\label{E: rearr}
    \begin{split}
        p>&\ppp^n[|X_{t\cj\sigma}^n|_A\leq r,\;\forall t\geq n]\\
        \geq  & \ppp^n[\tau ^n<n\;\text{and}\; |X_{t\cj\sigma}^n|_A\leq  r,\;\forall t\geq n]\\
        =& \ppp^n[\tau^n<n]\ppp[|X_{t\cj\sigma}^n|_A\leq r,\;\forall t\geq n\;|\;\tau^n<n]\\
        \geq& \ppp^n[\tau^n<n]\ppp[|X_{t\cj\sigma}^n|_A\leq r,\;\forall t\geq \tau^n],\\
    \end{split}
\end{equation}
By the definition of Pr-U.S. in view of Remark \ref{rem: epsdelta}, there exists an $\eta=\eta(r,\eps)<r$ such that
$\ppp^n[|X_{t\cj\sigma}^n|_A\leq r,\;\forall t\geq \tau^n]\geq \eps$. We choose $\eps$ sufficiently close to $1$ so that, by \eqref{E: rearr},
\begin{equation}\label{E: conf}
   \ppp^n[\tau^n<n]=p-\hat{p}<p, 
\end{equation}
where $\hat{p}=\hat{p}(\eps)\ll 1$. Note that we have implicitly defined $\eta$ and $\tau^n$ based on the choice of $\eps$ such that \eqref{E: conf} holds. 

However, by Remark \ref{rem: nonexplos} and Proposition \ref{thm: weak_compact}, there exists a subsequence, still denoted by 
$X^n\in\soln(x_n,W)$, such that $x_n\rightarrow x$ and $(X^n)^\sigma\rightharpoonup X^\sigma$ with $X\in\soln(x,W)$ on any compact interval of $\R_{\geq  0}$. By Skorohod \cite[Theorem 2.4]{da2014stochastic}, there exists a probability space $(\tilde{\Omega}^\dagger, \tilde{\mathscr{F}}^\dagger, \{\tilde{\mathscr{F}}_t^\dagger\},\tilde{\pp}^\dagger)$, a process $\tilde{X}^\sigma$  and a sequence of processes  $\{(\tilde{X}^n)^\sigma\}$ with laws $\ppp$ and $\{\ppp^n\}$, respectively, such that  
\begin{equation}\label{E: asconv}
    \lim\limits_{n\rightarrow\infty}(\tilde{X}^n)^\sigma=\tilde{X}^\sigma,\;\; \tilde{\pp}^\dagger-\text{a.s.}.
\end{equation}
Let $\tau=\inf\{t\geq 0: |\tilde{X}_{t\cj\sigma}|_A\leq \eta/2\}$,
due to the asymptotic behavior, we have
\begin{equation}\label{E: asymp2}
    \tilde{\pp}^\dagger[\tau<\infty]\geq p.
\end{equation}
By \eqref{E: asconv} and \eqref{E: asymp2}, there exists some sufficiently large $N_1(\eta,q_1)$ and $N_2(\eta,q_2)$  such that for any arbitrary $q_1,q_2\in(0,1)$,
\begin{equation}\label{E: conv}
    \tilde{\pp}^\dagger\left[\sup_{t\in[0,N_1]}|\tilde{X}^n_{t\cj\sigma}-\tilde{X}_{t\cj\sigma}|\leq\eta/2\right]\geq q_1,
\end{equation}
\begin{equation}\label{E: tau2}
    \tilde{\pp}^\dagger\left[\tau<N_2\right]\geq pq_2.
\end{equation}
Note that the events in  \eqref{E: conv} and \eqref{E: tau2} are independent, combining these and choosing $n\geq \max(N_1,N_2)$, we have
\begin{equation*}
    \begin{split}
       &\tilde{\pp}^\dagger\left[\exists t<n\;\text{s.t.}\;\tilde{X}^n_{t\cj\sigma}\in\oball_\eta(A)\right]
       \geq  \tilde{\pp}^\dagger\left[Q\right]
       \geq  pq_1q_2,
    \end{split}
\end{equation*}
where 
\begin{equation*}
    \begin{split}
        Q&:=\{\sup_{t\in[0,N_1]}|\tilde{X}^n_{t\cj\sigma}-\tilde{X}_{t\cj\sigma}|\leq\eta/2\}\\
        &\cap\{|\tilde{X}_{\tau\cj\sigma}|_A\leq \eta/2\;\text{for}\;\tau<N_2\}.
    \end{split}
\end{equation*}
We let $q_1q_2>\frac{p-\hat{p}}{p}$, then there exists an $n$ such that
\begin{equation}\label{E: conflict}
    \ppp^n[\tau_n<n]> p-\hat{p}.
\end{equation}
Equation \eqref{E: conflict} contradicts \eqref{E: conf}, which completes the proof.
\end{pf}

\begin{cor}
If $\sys$ satisfies a stability with safety guarantee specification $(\init,A, U, p)$ and $\init$ is compact, then for every $\eps>0$, $\sys$ satisfies the reach-avoid-stay specification $(\init, \ball_\eps(A), U, p)$.
\end{cor}
\begin{pf}
 We add the condition $\{\sigma=\infty\}$, then (1) and (2) are still equivalent in Proposition \ref{prop: uniform}. The conclusion follows directly by the definitions of the two specifications. 
\end{pf}

\subsection{\bf The converse side}
The converse side is intended to show probabilistic stability with safety is  necessary to probabilistic reach-avoid-stay specifications. Unfortunately, due to the diffusion effects and
the concept of weak solutions, probabilistic reach-avoid-stay
specifications, other than  reach-avoid-stay with probability
one, may fail to be related to probabilistic stability
with safety guarantees w.r.t. some subset of the target set. Considering our core target is to show the sufficiency,  we  convey the main idea and show the proofs in Appendix.  

Throughout this subsection, we 
suppose that $\init,\Gamma\subseteq\dd$ and $\sys$ satisfies a reach-avoid-stay specification $(\init,\Gamma, U, p)$. We first make a quick judgement that there exists a probability-$p$ invariant compact subset of $\Gamma$. 
\begin{lem}\label{lem: exist}
Suppose that $\Gamma$ is compact and $\init$ is nonempty. If $\sys$  satisfies a reach-avoid-stay specification $(\init, \Gamma, U, p)$ with $p\in(0,1]$, then the set
\begin{equation*}
    A=\{x\in\Gamma: \forall X\in\soln(x,W), \ppp^X[X_{t}\in\Gamma,\;\forall t\geq 0]\geq p\}
\end{equation*}
is a nonempty and compact set  with 
$\ppp^X[X_{t}\in A,\;\forall t\geq 0]\geq p$ for all $X\in\soln(A,W)$.
\end{lem}

The next lemma shows that given an arbitrary solution $X$ of $\sysp$, we can construct a weak solutions $Z$ for $\sys$ that solves the martingale problem and is relatively close to $X$. 

\begin{lem}\label{lem: construct}
Let $\delta'\in(0,\delta)$ 
and  
 $\tau$ be such that $\tau<\sigma$ a.s.. 
  Then 
  there exists some $r=r(\tau, \delta, \delta')$ such that for every $X\in\solnp(x,W)$ with $x\in \mathcal{D}$ and for all $z\in\ball_r(x)$, there exists a weak solution $Z\in\soln(z,W)$
  such that 
  $Z_{T\cj\sigma}\in\ball_r(X_{T\cj\sigma})$ a.s. for $T\in[\tau,\infty)$.

\end{lem}

It can be shown that under the construction of Lemma \ref{lem: exist} and  \ref{lem: construct}, the set $A$ (generated by solutions of $\sys$) is Pr-U.A. for any weak solution of $\sysp$ with $\delta'\in (0,\delta)$. 
\begin{prop}\label{prop: attract}
Suppose that $\sys$ satisfies reach-avoid-stay specification $(\init, \Gamma, U, p)$. 
Let $A$ be the set given in Lemma \ref{lem: exist}. Then $A$ is uniformly attractive for $\sysp$ with probability at least $p$, i.e., for every $\eps>0$, there exists $T=T(\eps,p)$ such that for any $X\in\bigcup_{x\in \init}\solnp(x,W)$,
    $$\ppp\left[|X_{t\cj\sigma}|_A<\eps,\;\forall t\geq T\right]\geq p. $$
\end{prop}

However, for non-strictly invariant sets ($p<1$), we are not able to show the Pr-U.S. property due to a geometric gap where we cannot  arbitrarily set $\eps$ and $r$ as in Definition \ref{def: u.s.}. On the other hand, if there exists a strictly invariant subset of $\Gamma$, nice properties appear. This is not a surprise given  \cite[Proposition 16]{liu2020smooth}. The possibility of such an existence occurs when the system admits a family of a.s. stable Dirac invariant measures for each signal $\xi$, which are strictly contained in $\Gamma$. We convert the statement into the stochastic context in the next proposition. The proof 
relies on a similar construction  \cite[Lemma 15]{liu2020smooth} as Lemma \ref{lem: construct}. We skip the proofs for the following two statements.

\begin{prop}
For system $\sys$,  any nonempty compact set $A\subseteq \dd$ with 
$\ppp^X[X_{t}\in A,\;\forall t\geq 0]=1$ is Pr-U.A.S. for $\sysp$ whenever $\delta'\in[0,\delta)$.
\end{prop}
\begin{cor}
If $\sys$ satisfies a reach-avoid-stay specification $(\init, \Gamma, U, 1)$ with compact $\init$, then there exists a nonempty compact set $A\subseteq \Gamma$ with 
$\ppp^X[X_{t}\in A,\;\forall t\geq 0]=1$ such that for  any $\delta'\in(0,\delta)$, $\sysp$ satisfies a stability with safety specification $(\init, A, U, 1)$.
\end{cor}

\section{Lyapunov-Barrier Conditions for Probabilistic Stability with Safety}\label{sec: lb}
We aim to show how Lyapunov-Barrier functions can sufficiently guarantee the probabilistic stability with safety in this subsection. 
Recall region $\dd$ and $\rr$ in Definition \ref{def:workspace}.
\begin{deff}[Stochastic Lyapunov functions] 
Let $A\subset\dd$ be a closed set. A function $V\in (C^2(\aba);\R_{\geq 0})$ is said to be a stochastic Lyapunov functions (SLF) w.r.t. $A$ if there exist $\alpha_1,\alpha_2,\alpha_3\in\mathcal{K}$ such that, for all $x\in \aba$, 
   \begin{equation}\label{E: V1}
    \alpha_1(|x|_A)\leq V(x)\leq \alpha_2(|x|_A)
\end{equation}
and 
\begin{equation}\label{E: V2}
    \sup\limits_{d\in\ball_\delta}\LL_d V(x)\leq -\alpha_3(|x|_A).
\end{equation}

\end{deff}

We first make a quick extension of the existing Lyapunov theorems to systems with point mass perturbations.

\begin{lem}[Uniform recurrence]\label{lem: reccur}
Given an SLF $V$, there exists some $\eta>0$ such that, for every $\eps\in(0,1)$ and $r\in (0,R/2)$, there exists some $T=T(\eps,\eta,  r)>0$ such that, for any $x\in\oball_\eta(A)$, we have
$\infx\ppp^X[\tau<T ]\geq 1-\eps,$
where $\tau=\inf\{t\geq 0: X_t\in \oball_r(A)\}$ is the first hitting time of $\oball_r(A)$ for each $X\in\soln(x,W)$.
\end{lem}

The proof of Lemma \ref{lem: reccur} is completed in the Appendix. 
\begin{prop}
Suppose $A\subset \dd$ is compact. If there exists an SLF $V$ w.r.t. $A$, then $A$ is Pr-U.A.S for $\sys$.
\end{prop}
\begin{pf}
By a standard supermartingale argument \cite[Lemma 1, Chap II]{kushner1967stochastic}, we can show that the existence of SLF  implies Pr-U.S.. 
To show Pr-U.A., let $r\in(0,R/2)$, then by Pr-U.S. and Remark \ref{rem: epsdelta}, there exists a $\kappa\in(0,r)$ such that 
$\supx\ppp^X[|X_{t}|_A\leq r\;\forall t\geq s(\omega)]\geq 1-\frac{\eps}{2}$
whenever $|X_s|_A\leq \kappa$. Now let $\tau=\inf\{t\geq 0: X_{t}\in \oball_\kappa(A)\}$. By Lemma \ref{lem: reccur}, there exists some $\eta>0$ such that  we can find a $T=T(\eps/2, \eta, \kappa)$ to make 
$\infx\ppp^X[\tau<T]\geq 1-\frac{\eps}{2}. $
Therefore, for all $|x|_A<\eta$ and for all $X\in\soln(x,W)$,
\begin{equation}
    \begin{split}
     &\ppp^X[|X_t|_A\leq r,\;\forall t\geq T]\\
        \geq  & \ppp^X[\tau <T\;\text{and}\; |X_t|_A\leq  r,\;\forall t\geq T]\\
        \geq& \ppp^X[\tau<T]\ppp^X[|X_t|_A\leq r,\;\forall t\geq \tau\;|\;\tau<T]\\
        \geq& \ppp^X[\tau<T](1-\eps/2)
        \geq  (1-\eps/2)^2\geq 1-\eps.
    \end{split}
\end{equation}
\end{pf}

The following result demonstrates that the existence of  SLFs is sufficient to characterize probabilistic stability with safety  specifications with probabilities depending on initial conditions. 

\begin{thm}\label{thm: Lya}
Suppose that $A\subset\dd$ is compact and $\aba\subset\rr$. If there exists an SLF $V\in(C^2(\aba);\R_{\geq 0})$ and some $G:=\oball_r(A)$ such that (i) $r\in(0,R]$ and  $G\subset\dd$; (ii) $\init\subset G$,
then $\sys$ satisfies the probabilistic stability with safety  specification $(\init,A, U, 1-\frac{\sup_{x\in\init}V(x)}{\alpha_1(r)})$.
\end{thm}

We need the following lemma to accomplish the proof. 

\begin{lem}\label{lem: asymp}
For each $X\in\bigcup_{x\in\init}\soln(x,W)$, set $\tau:=\inf\{t\geq 0: X_t\in G^c\}$. Then for all $X\in\bigcup_{x\in\init}\soln(x,W)$, 
$$\ppp^X[\lim\limits_{t\rightarrow\infty}|X_t|_A=0\;|\;\tau=\infty]=1.$$
\end{lem}
\begin{pf}
By a similar approach to Lemma \ref{lem: reccur}, we set arbitrary $r^*\in (0,r)$. By the Pr-U.S. property, for all $X\in\soln(x,W)$ there should exist $\eta\in(0,r^*)$ such that for any $\eps\in(0,1)$,  $X_{\tau^*}\in\ball_\eta(A)$ implies $\ppp^X[X_t\in\oball_{r^*}(A),\;\forall t\geq \tau^*]\geq 1-\eps$, where 
$
    \tau^*=\inf\{t\geq 0: X_t\in \oball_{\eta}(A)\}.
$
By It\^{o}'s formula, for each weak solution  we have
\begin{equation}\label{E: flowout}
    \begin{split}
        0 &\leq V(x)+\eee^X\int_0^{\tau^*\cj\tau\cj t}\LL_d V(X(s))ds\\
        &\leq V(x)-\alpha_3(\eta)\eee^X[\tau^*\cj\tau\cj t]\\
    \end{split}
\end{equation}
Since that on $\{\tau^*\cj\tau\geq t\}$ we have $\tau^*\cj\tau\cj t=t$, thus
\begin{equation*}
    \begin{split}
        \eee^X[\tau^*\cj\tau\cj t]& \geq \int_\Omega \mathds{1}_{\{\tau^*\cj\tau\geq t\}}\cdot t\; d\ppp^X(\omega)=t\ppp^X[\tau^*\cj\tau\geq t],
    \end{split}
\end{equation*}
combining with \eqref{E: flowout} we have
\begin{equation}
    \ppp^X[\tau^*\cj\tau\geq t]\leq V(x)/t\alpha_3(\eta),\;\;\text{for each}\; t,
\end{equation}
which implies $\ppp^X[\tau^*\cj\tau<\infty]=1$ for all $X\in\bigcup_{x\in\init}\soln(x,W)$. On $\{\tau=\infty\}$, for all weak solution, we have $\ppp^X[\tau^*<\infty]=1$  and 
\begin{align*}
&\ppp^X[\limsup_{t\rightarrow\infty}|X_t|_A\leq r^*]\\
\geq &\ppp^X[|X_t|_A\leq r^*,\;\forall t\geq \tau^*\;|\;\tau^*<\infty]\geq 1-\eps. 
\end{align*}
Since $\eps$ and $r^*$ are arbitrary, the conclusion follows. 
\end{pf}
\begin{rem}
Lemma \ref{lem: asymp} shows that SLFs eliminate the possibility of safe sample paths up/down-crossing any neighborhood of $A$ infinitely often. \cite[Theorem 2, Chap II]{kushner1967stochastic} demonstrates the same result by constructing the total time spent in $G\setminus\oball_\eps(A)$ after time $t$ and showing that it converges  a.s. to $0$ as $t\rightarrow\infty$.
\end{rem}

\textbf{Proof of Theorem \ref{thm: Lya}}.
The existence of SLF shows that $A$ is Pr-U.A.S. for $\sys$. Now, for all $X\in\soln(x,W)$ with $x\in\init$, define  $\tau:=\inf\{t\geq 0: X_t\in G^c\}$.
Then, for all $t\geq 0$ and for all $X\in\soln(x,W)$, \begin{equation}\label{E: leq}
    \eee^X [V(X_{\tau\wedge t})]=V(X_0)+\eee^X\left[\int_0^{\tau\cj t} \LL_d V(X_s)ds\right]\leq V(x),
\end{equation}
and, for all $t\geq 0$, 
\begin{equation}\label{E: geq}
    \eee^X [V(X_{\tau\cj t})]\geq \eee^X [\mathds{1}_{\{\tau\leq t\}}V(X_{\tau})]> \alpha_1(r)\ppp^X[\tau\leq t],
\end{equation}
which imply
\begin{equation}
    \ppp^X[\tau\leq t]< V(x)/\alpha_1(r),\;\;\forall t\geq 0.
\end{equation}
Sending $t\rightarrow\infty$ we get for all $X\in\bigcup_{x\in \init}\soln(x,W)$, 
$\ppp^X[\tau<\infty]<V(x)/\alpha_1(r)$, i.e.,
\begin{equation}
  \inf\limits_{X\in\soln(x,W)}\ppp^X[\tau=\infty]\geq  1-\frac{\sup_{x\in\init}V(x)}{\alpha_1(r)}.  
\end{equation}
Since $\{\tau=\infty\}\subset\{\sigma=\infty\}$ and by Lemma \ref{lem: asymp}, the conclusion follows. 
\pfbox

We have seen in the proof that conditions $\alpha_1(x)\leq V(x)$ and $\sup_{d\in\ball_\delta}\LL_d V\leq 0$ play the role of guaranteeing the probabilistic set invariance. We refer these conditions as the stochastic barrier certificates. An application in control synthesis, termed as stochastic control barrier functions,  has been shown in \cite[Proposition III.8]{wang2021barrier} with better safety probability compared to the zeroing-type barrier certificates \cite{prajna2007framework}
, however, less effectiveness than the reciprocal-type barrier certificates (see the definition in \cite[Definition III.1]{wang2021barrier}). A thorough comparison between the above mentioned stochastic barrier functions can be found in \cite{wang2021barrier}.
To provide stability with safety with a higher probability, one can combine SLF with the reciprocal-type barrier functions.

\begin{thm}\label{thm: lb}
Under the same assumption in Theorem \ref{thm: Lya}. Suppose there exists an SLF $V\in(C^2(\aba);\R_{\geq 0})$, some $G:=\oball_r(A)$ such that $G\in\dd$ and $\init\subset G$, as well as a function $B\in (C^2(G);\R_{\geq 0})$ satisfying 
\begin{enumerate}
    \item[(i)]  $\exists \alpha_1,\alpha_2\in\mathcal{K}$ s.t.
    \begin{equation}\label{E: b1}
        \frac{1}{\alpha_1(|x|_A)}\leq B(x)\leq \frac{1}{\alpha_2(|x|_A)}, \;\;\forall x\in G;
    \end{equation}
   
    \item[(ii)] $\exists \alpha_3\in\mathcal{K}$ s.t.
    \begin{equation}\label{E: b2}
        \sup\limits_{d\in \ball_\delta}[\LL_d B(x)-\alpha_3(|x|_A)]\leq 0, \;\;\forall x\in G.
    \end{equation}
   
\end{enumerate}
Then $\sys$ satisfies the probabilistic stability with safety  specification $(\init,A, U, 1)$.
\end{thm}
\begin{pf}
The proof is similar to Theorem \ref{thm: Lya}. We rely on the SLF to provide the property shown in Lemma \ref{lem: asymp}. Then the reciprocal type barrier function $B$ guarantees that $\ppp[\tau=\infty]=1$ \cite[Theorem 1]{clark2019control} for each weak solution.
\end{pf}
\begin{rem}\label{rem: relax}
Suppose $U^c:=\{x\in\R^n: h(x)\geq 0\}$ where $h$ is smooth, one can possibly enlarge $G$ such that $G\cap U\neq \emptyset$ with $\partial (G\cap U) $ being piecewise smooth. To see the satisfaction of  stability with safety specifications, along with the old conditions, one can introduce an extra reciprocal-type barrier function, denoted by $\tilde{B}$, and verify extra conditions that are similar to \eqref{E: b1} and \eqref{E: b2} by replacing $|x|_A$ with $h(x)$. 
\end{rem}

\section{Applications in Control Problems}\label{sec: app}

In this section, based on the results from Section \ref{sec: conn} and \ref{sec: lb}, we make a straightforward extension to a stochastic control Lyapunov-barrier  characterization for $\sys$ satisfying a probabilistic reach-avoid-stay specification $(\init,\Gamma, U, p)$ under controls. As a continuation of \cite{meng2021control}, we conduct a case study on enhancing the performance of jet engine compressors, under both noisy disturbances and bounded point mass perturbations, based on a reduced Moore-Greitzer nonlinear SDE model. 

\subsection{\bf Probabilistic reach-avoid-stay control via stochastic control Lyapunov-barrier functions}
We first recast the notion from Section \ref{sec: pre} for control systems. Given a nonempty compact convex set of control inputs  $\uu\subset\mathbb{R}^p$, 
consider a nonlinear system of the form
\begin{equation}\label{E: sys2}
dX_t=f(X_t)dt+\xi(t)dt+b(X_t)udt+g(X_t)dW_t, \;X_0=x,
\end{equation}
where the mapping $b:\mathbb{R}^n\rightarrow\R^{n\times p}$ is smooth;  $u:\R_{\geq 0}\rightarrow \uu$ is a locally bounded measurable control signal, whilst the other notation remains the same.

\begin{deff}[Control strategy]
A control strategy is a function
\begin{equation}
   \kappa:\R^n\rightarrow \uu. 
\end{equation}

\end{deff}
We further denote $\sys^\kappa$ by the control system driven by \eqref{E: sys2} that is comprised by $u=\kappa(x)$.

\begin{deff}
(\textbf{Probabilistic reach-avoid-stay controllable}): Given $\init,\Gamma\subseteq \dd$ 
and $p\in [0,1]$,  $\sys$ is said to be probabilistic  reach-avoid-stay controllable w.r.t.  $(\init,\Gamma, U, p)$, if there exists a Lipschitz continuous control strategy $\kappa$ such that the system $\sys^\kappa$ satisfies the specification $(\init,\Gamma, U, p)$. 
\end{deff}
\begin{prop}\label{prop: control}
Given $\init,\Gamma\subseteq \dd$, if there exists a smooth function $V\in(C^2(\aba);\R_{\geq 0})$ and some  $G:=\oball_r(A)$, such that
\begin{enumerate}
    \item[(i)]
    $r\in(0,R]$, $G\subset \dd$ and $\init\subset G$;
    \item[(ii)] 
    $
    \alpha_1(|x|_A)\leq V(x)\leq \alpha_2(|x|_A)
$
and
    $$
    \inf\limits_{u\in \uu}\sup\limits_{x\in S}\sup\limits_{d\in\delta\mathcal{B}}[\LL_d^u V(x)+\alpha_3(|x|_A)]\leq  0,
$$
for some $\alpha_1,\alpha_2,\alpha_3\in\mathcal{K}$, where 
$\LL_d^uV(x):=\LL_d V(x)+\nabla V(x)\cdot b(x)u.$
\end{enumerate}
Then $\sys$ is probabilistically reach-avoid-stay controllable  w.r.t. $(\init,\Gamma, U,1-\frac{\sup_{x\in\init}V(x)}{\alpha_1(r)})$. 
\end{prop}

Similarly, one can extend the above proposition to find sufficient conditions for a `probability 1' reach-avoid-stay  based on Theorem \ref{thm: lb}. Apart from the conditions in Proposition \ref{prop: control}, one need to additionally verify if there exists a $B\in (C^2(G);\R_{\geq 0})$ satisfying 
\begin{enumerate}
    \item[(i)] 
    $ \frac{1}{\tilde{\alpha}_1(|x|_A)}\leq B(x)\leq \frac{1}{\tilde{\alpha}_2(|x|_A)}, \;\;\forall x\in G
$
    for some  class-$\mathcal{K}$ functions $\tilde{\alpha}_1,\tilde{\alpha}_2$; 
    \item[(ii)] 
$
    \inf\limits_{u\in \uu}\sup\limits_{x\in S}\sup\limits_{d\in\delta\mathcal{B}}[\LL_d^u B(x)-\tilde{\alpha}_3(|x|_A)]\leq  0,
$
   for some class-$\mathcal{K}$ function $\tilde{\alpha}_3$.
\end{enumerate}

\begin{rem}
In view of Remark \ref{rem: relax}, the region $G$ can be further relaxed if $\partial U$ is smooth enough. Correspondingly, some extra conditions given by another reciprocal control barrier function are needed to guarantee the sufficiency of $(\init,\Gamma,U,1)$ controllability.
\end{rem}

\subsection{\bf Case study}
We use the reduced Moore-Greitzer SDE model with an additive control input $[v,0]^T$ and a multiplicative noise to illustrate the effectiveness. The model is given as:
\begin{equation}\label{E: MG}
\begin{split}
 \dfrac{d}{dt}\begin{bmatrix}
\Phi(t)\\ \Psi(t)
\end{bmatrix}& =\begin{bmatrix}
\frac{1}{l_c}(\psi_c-\Psi(t))\\
\frac{1}{16l_c}\left(\Phi(t)-\mu\sqrt{\Psi(t)}\right)
\end{bmatrix}+\begin{bmatrix}
\xi_1(t)\\\xi_2(t)
\end{bmatrix}\\
& +\eps\begin{bmatrix}
(\Phi(t)-\Phi_e(\mu))\beta_1(t) \\  (\Psi(t)-\Psi_e(\mu))\beta_2(t)
\end{bmatrix}+\begin{bmatrix}
v(t)\\0
\end{bmatrix},
\end{split}
\end{equation}
where $\psi_c=a+\iota*[1+\frac{3}{2}(\frac{\Phi}{\Theta}-1)-\frac{1}{2}(\frac{\Phi}{\Theta}-1)^3]$, $\beta_1,\beta_2$ are i.i.d. Brownian motions, $(\Phi_e(\mu),\Psi_e(\mu))=:X_e(\mu)$ are equilibrium points for $\xi_1,\xi_2,v\equiv 0$. The other parameters are as follows:
\begin{equation*}
    l_c=8,\;\iota=0.18,\;\Theta=0.25,\;a=0.67\iota,\;\eps=0.08,\;\delta=0.01.
\end{equation*}
The physical meanings of variables, parameters and the description of this model can be found in \cite[Section V]{meng2021control}.\\
\begin{rem}
For $\xi_1,\xi_2,v\equiv 0$, the system admits a family of equilibrium points $X_e(\mu)$ depending on the tunable parameter $\mu$. As $\mu$ drops in the neighborhood of the deterministic Hopf bifurcation point, the system undergoes a D-bifurcation (the  stability of the invariant measure $\delta_{\{X_e\}}$ changes and a new invariant measure in $\R^n\setminus\{X_e\}$ is built up) and a P-bifurcation (the shape of density of the new measure changes). 
The full stochastic Hopf bifurcation diagram in \cite[Fig 9.13]{arnold1995random} conveys the brief idea. 

\begin{figure}[h]
\centering
\includegraphics[scale=0.27]{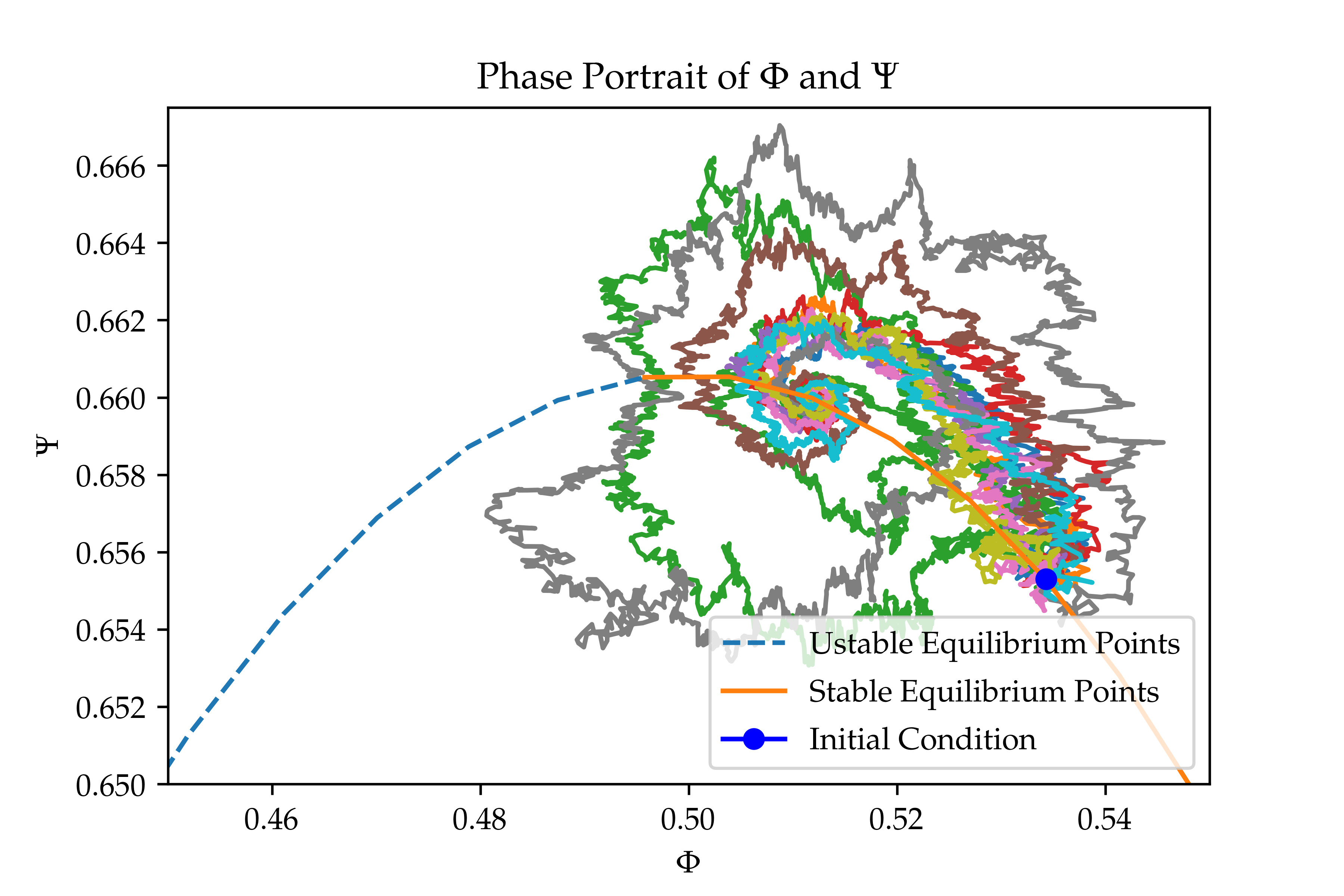}
\includegraphics[scale=0.27]{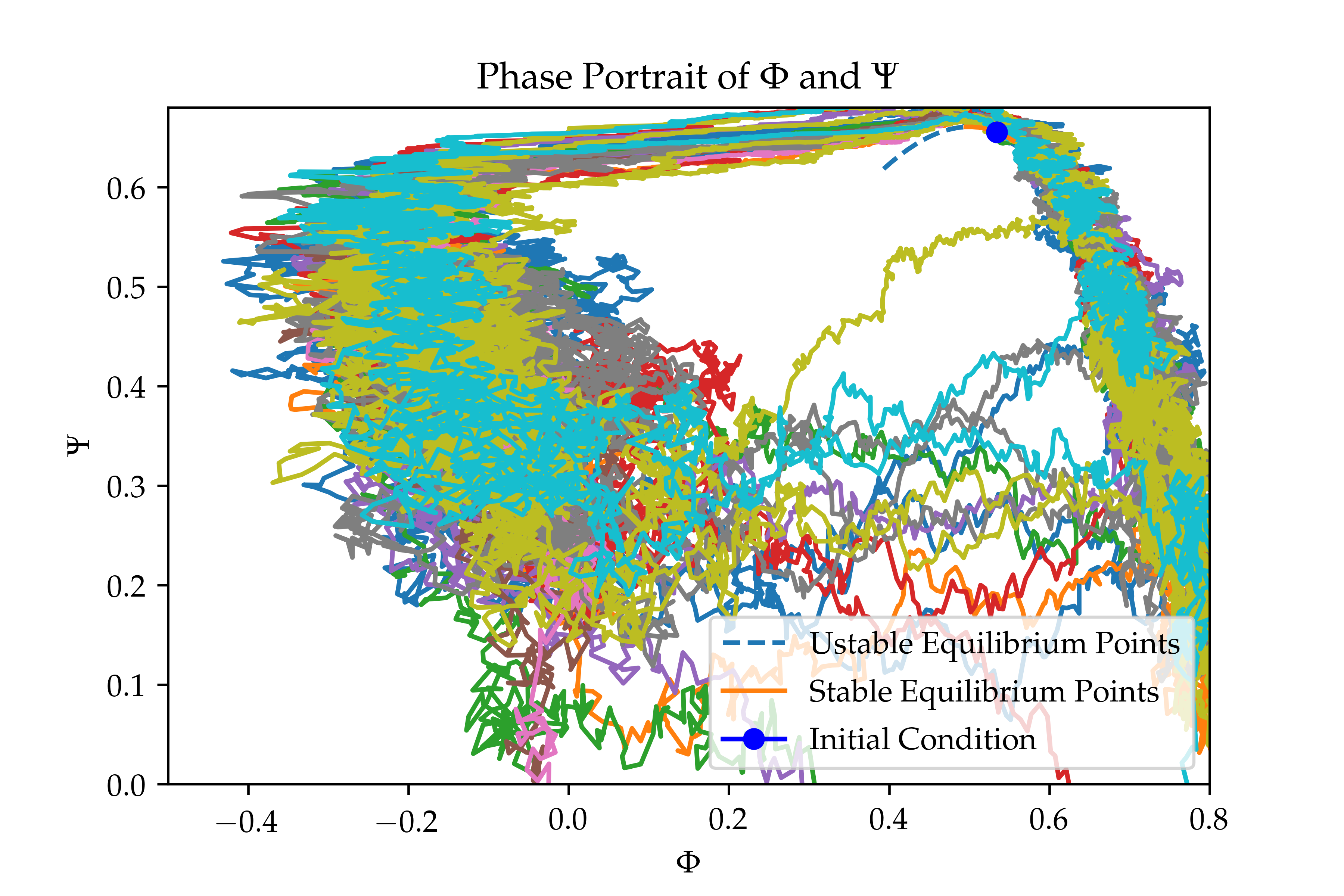}
\caption{Given the initial condition of $(\Phi_0,\Psi_0)=(0.5343,0.6553)$, the phase portraits are generated by \eqref{E: MG} with $\mu=0.63$ (left) and $\mu=0.59$ (right) under $\xi_1,\xi_2,v\equiv 0$.
}\label{F:phase_portrait}
\end{figure}

\begin{figure}[h]
\centering
\includegraphics[scale=0.27]{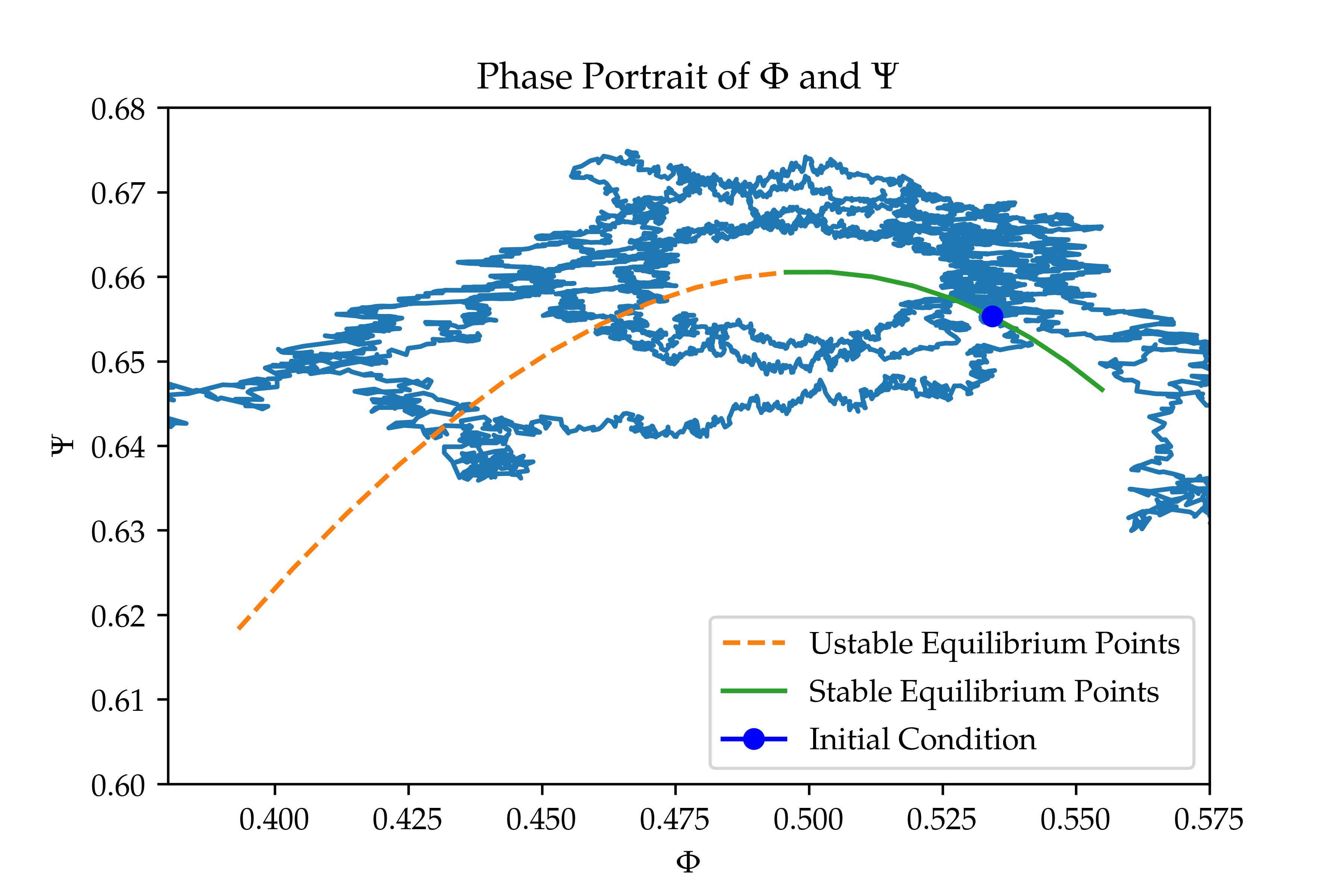}
\includegraphics[scale=0.27]{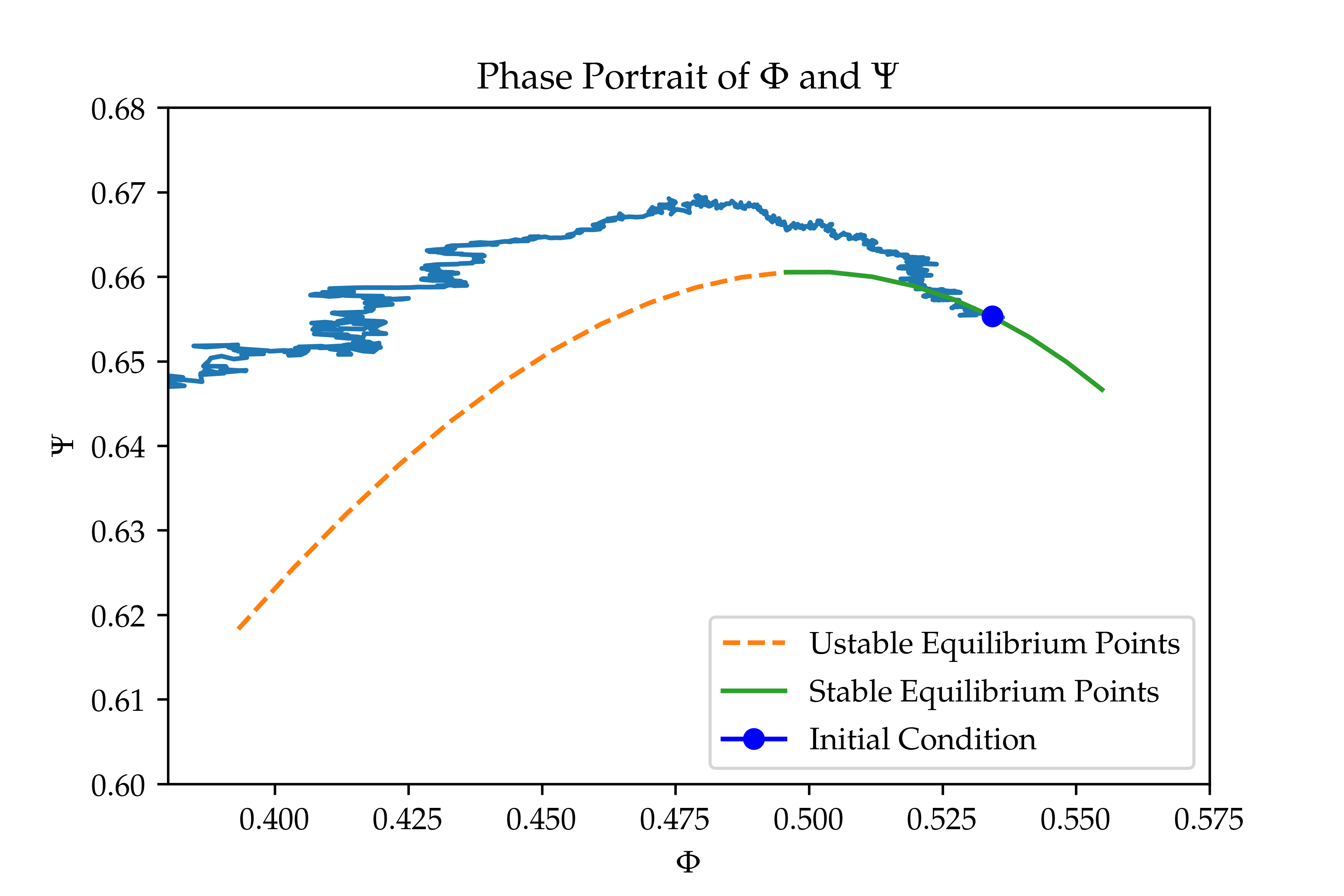}
\caption{For $(\Phi_0,\Psi_0)=(0.5343,0.6553)$, two sample paths (left: converges; right: diverges) are generated with $\mu=0.609$ (in between the $1^{\text{st}}$ D-bifurcation point and P-bifurcation point) under $\xi_1,\xi_2,v\equiv 0$.
}\label{F:phase_portrait2}
\end{figure}

Within the a.s. exponentially stable region, any bounded perturbation $\xi$ causes a bounded long-term perturbation of $X_e(\mu)$, and ultimately formulate a compact set containing $X_e(\mu)$. For unstable $\delta_{\{X_e\}}$, especially for those after P-bifurcation, we are interested in stabilizing the robust system to a compact set. 
\end{rem}

\begin{prob}\label{P: case_study}
We aim to manipulate  $\mu$ and $v$ simultaneously such that the state $(\Phi,\Psi)$ are regulated to satisfy reach-avoid-stay specification $(\init,\Gamma,U,1)$. We require that $\mu:\R_{\geq 0}\rightarrow [0.5,1]$ is time-varied with $\mu(0)\in[0.62,0.66]$ and  $|\mu(t+\tau)-\mu(t)|\leq 0.01\tau$ for any $\tau>0$.  We define $\init=\{(\Phi_e(\mu(0)),\Psi_e(\mu(0)))\}$; $\Gamma$ to be the ball that centered at $\gamma=(0.4519,0.6513)$ with radius $r=0.013$, i.e. $\Gamma=\gamma+r\ball$; the unsafe set $U=\{(x,y):h_1\leq 0\}\cap \{(x,y):h_2\leq 0\}$, where $h_1(x,y)=-|(x,y)-(0.49,0.64)|+0.055$, $h_2(x,y)=|(x,y)-(0.50,0.65)|-0.003$. We set 
$v\in\uu=[-0.05,0.05]\cap \R$.
\end{prob}

We refer readers to \cite[Remark 30]{meng2021control} for treatments of $\mu$ as another control input. For each SDE, the signals $\xi_1,\xi_2$ of each sampling time is generated randomly from $\{-0.1,0.1\}$.

We choose SLF $V(x,y)=\frac{l_c}{2}(x-\gamma_1)^2+8l_c(y-\gamma_2)^2$ and $\alpha_3(x)=0.1x$; set $B_{i}=-\log\left(\frac{h_i}{1+h_i}\right)$ for $i=1,2.$ The settings for the quadratic programming keep the same as \cite[Section V.B]{meng2021control}. We mix sample paths under different $\xi_1,\xi_2$  and show the simulation results as below:

\begin{figure}[h]
\centering
\includegraphics[scale=0.27]{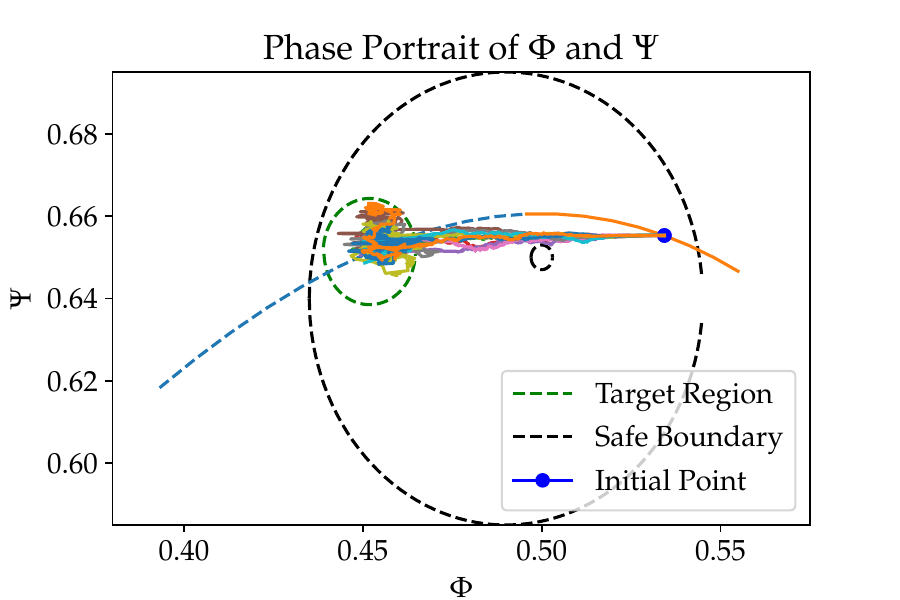}
\includegraphics[scale=0.27]{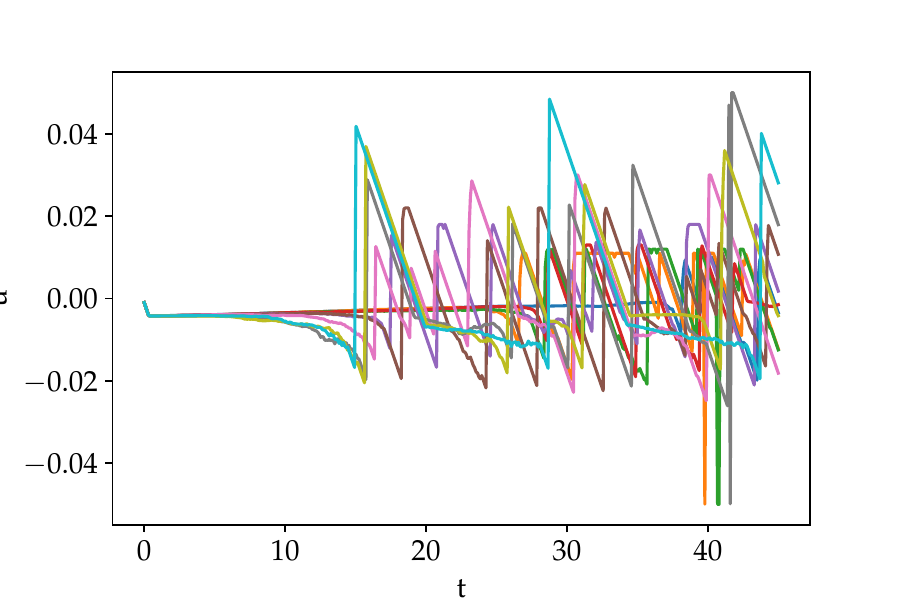}
\includegraphics[scale=0.27]{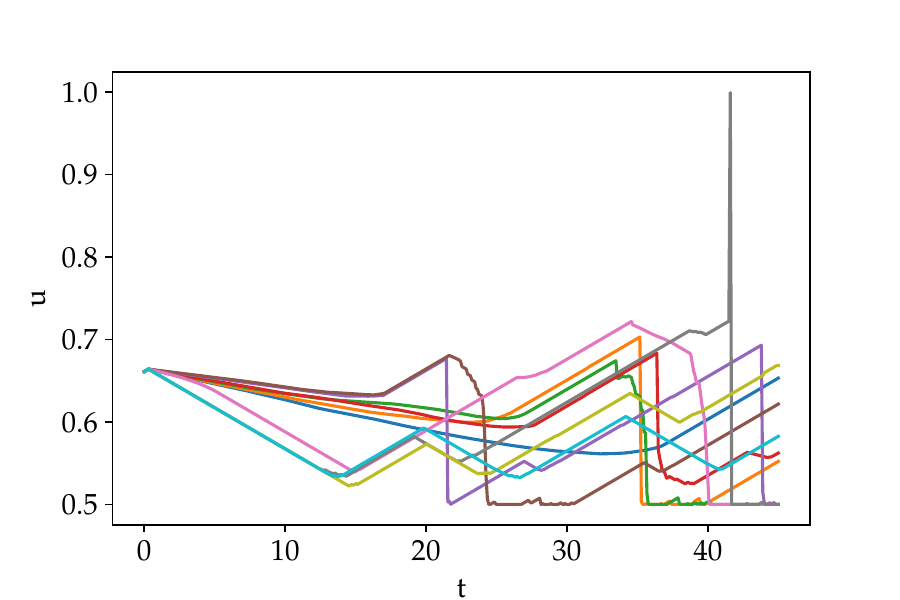}
\caption{Controlled sample paths, control input $v,\mu$ of  Problem \ref{P: case_study}.
}\label{F:num1}
\end{figure}
\begin{rem}
Note that we have adopted reciprocal type barrier functions, which potentially generates impulse-like control signals (to cancel the diffusion effects) and terminates the programming. However, once the synthesis succeeds, the feasible controlled sample paths satisfy the specification.
\end{rem}

\section{CONCLUSIONS}\label{sec: conclusion}
In this paper, we formulated stochastic Lyapunov-barrier functions to develop sufficient conditions on probabilistic reach-avoid-stay specifications. Given uncertainties of the model, robustness was taken into account such that a worst-case scenario is guaranteed.  We characterized a general topological structure of the initial sets, target sets and unsafe sets under the stochastic settings and discussed relaxations given the smoothness of the unsafe boundary.  We investigated the effectiveness in a case study of jet engine compressor control problem. Despite of the potentially unbounded control inputs, the control version of SLF along with reciprocal-type barrier functions guarantee a probability-$1$ satisfaction. 

However, just like deterministic  Lyapunov-like functions only providing a stability characterization of the solutions, the stochastic Lyapunov-type argument can only estimate a lower bound of `satisfaction in probability/law' without solving the evolving states and distributions. It renders more difficulties of selecting Lyapunov/barrier functions under the restrictive geometric requirements of the initial conditions and unsafe sets. 

For future work, compared to the rough estimation of `probabilistic domain of satisfactions' given Lyapunov-like functions, 
it would be necessary to consider accurate evaluation of laws and  investigate formal methods in providing more reliable schemes on finding the probabilistic winning sets (from which the specifications are satisfied). Considering the uncertainties of stochastic modelling, to provide soundness and (possibly weak) completeness, stochastic abstraction analysis is fundamental to the robust stochastic control synthesis problems.

\bibliographystyle{plain}        
\bibliography{Manuscript}

@inproceedings{agrawal2017discrete,
  title={Discrete Control Barrier Functions for Safety-Critical Control of Discrete Systems with Application to Bipedal Robot Navigation.},
  author={Agrawal, Ayush and Sreenath, Koushil},
  booktitle={Proc. of RSS},
  year={2017}
}


@inproceedings{hsu2015control,
  title={Control barrier function based quadratic programs with application to bipedal robotic walking},
  author={Hsu, Shao-Chen and Xu, Xiangru and Ames, Aaron D},
  booktitle={Proc. of ACC},
  pages={4542--4548},
  year={2015},
  organization={IEEE}
}

@incollection{nguyen2020dynamic,
  title={Dynamic walking on stepping stones with gait library and control barrier functions},
  author={Nguyen, Quan and Da, Xingye and Grizzle, JW and Sreenath, Koushil},
  booktitle={Proc. of WAFR},
  pages={384--399},
  year={2020},
  publisher={Springer}
}

@article{nilsson2017augmented,
  title={Augmented finite transition systems as abstractions for control synthesis},
  author={Nilsson, Petter and Ozay, Necmiye and Liu, Jun},
  journal={Discrete Event Dynamic Systems},
  volume={27},
  number={2},
  pages={301--340},
  year={2017},
  publisher={Springer}
}

@inproceedings{faulwasser2009model,
  title={Model predictive path-following for constrained nonlinear systems},
  author={Faulwasser, Timm and Kern, Benjamin and Findeisen, Rolf},
  booktitle={Proc. of CDC},
  pages={8642--8647},
  year={2009},
  organization={IEEE}
}

@book{fribourg2013control,
  title={Control of switching systems by invariance analysis: applcation to power electronics},
  author={Fribourg, Laurent and Soulat, Romain},
  year={2013},
  publisher={John Wiley \& Sons}
}

@article{kloetzer2008fully,
  title={A fully automated framework for control of linear systems from temporal logic specifications},
  author={Kloetzer, Marius and Belta, Calin},
  journal={IEEE Transactions on Automatic Control},
  volume={53},
  number={1},
  pages={287--297},
  year={2008},
  publisher={IEEE}
}

@article{pola2008approximately,
  title={Approximately bisimilar symbolic models for nonlinear control systems},
  author={Pola, Giordano and Girard, Antoine and Tabuada, Paulo},
  journal={Automatica},
  volume={44},
  number={10},
  pages={2508--2516},
  year={2008},
  publisher={Elsevier}
}

@article{girard2009approximately,
  title={Approximately bisimilar symbolic models for incrementally stable switched systems},
  author={Girard, Antoine and Pola, Giordano and Tabuada, Paulo},
  journal={IEEE Transactions on Automatic Control},
  volume={55},
  number={1},
  pages={116--126},
  year={2009},
  publisher={IEEE}
}

@article{tabuada2006linear,
  title={Linear time logic control of discrete-time linear systems},
  author={Tabuada, Paulo and Pappas, George J},
  journal={IEEE Transactions on Automatic Control},
  volume={51},
  number={12},
  pages={1862--1877},
  year={2006},
  publisher={IEEE}
}

@inproceedings{liu2017robust,
  title={Robust abstractions for control synthesis: Completeness via robustness for linear-time properties},
  author={Liu, Jun},
  booktitle={Proc. of HSCC},
  pages={101--110},
  year={2017}
}


@book{baier2008principles,
  title={Principles of Model Checking},
  author={Baier, Christel and Katoen, Joost-Pieter},
  year={2008},
  publisher={MIT press}
}


@article{li2020robustly,
  title={Robustly complete synthesis of memoryless controllers for nonlinear systems with reach-and-stay specifications},
  author={Li, Yinan and Liu, Jun},
  journal={IEEE Transactions on Automatic Control},
  year={2020},
  publisher={IEEE}
}

@article{ames2016control,
  title={Control barrier function based quadratic programs for safety critical systems},
  author={Ames, Aaron D and Xu, Xiangru and Grizzle, Jessy W and Tabuada, Paulo},
  journal={IEEE Transactions on Automatic Control},
  volume={62},
  number={8},
  pages={3861--3876},
  year={2016},
  publisher={IEEE}
}

@article{romdlony2016stabilization,
  title={Stabilization with guaranteed safety using control Lyapunov--barrier function},
  author={Romdlony, Muhammad Zakiyullah and Jayawardhana, Bayu},
  journal={Automatica},
  volume={66},
  pages={39--47},
  year={2016},
  publisher={Elsevier}
}

@article{teel2000smooth,
  title={A smooth Lyapunov function from a class-$\mathcal{KL}$ estimate involving two positive semidefinite functions},
  author={Teel, Andrew R and Praly, Laurent},
  journal={ESAIM: Control, Optimisation and Calculus of Variations},
  volume={5},
  pages={313--367},
  year={2000},
  publisher={EDP Sciences}
}

@article{liu2020converse,
  title={Converse Barrier Functions via {L}yapunov Functions},
  author={Liu, Jun},
  journal={IEEE Transactions on Automatic Control},
  year={2021},
  publisher={IEEE}
}

@article{lin1996smooth,
  title={A smooth converse Lyapunov theorem for robust stability},
  author={Lin, Yuandan and Sontag, Eduardo D and Wang, Yuan},
  journal={SIAM Journal on Control and Optimization},
  volume={34},
  number={1},
  pages={124--160},
  year={1996},
  publisher={SIAM}
}

@article{sontag1989universal,
  title={A ‘universal’ construction of {A}rtstein's theorem on nonlinear stabilization},
  author={Sontag, Eduardo D},
  journal={Systems \& Control Letters},
  volume={13},
  number={2},
  pages={117--123},
  year={1989},
  publisher={Elsevier}
}

@article{braun2017existence,
  title={On (the existence of) control Lyapunov barrier functions},
  author={Braun, Philipp and Kellett, Christopher M},
  year={2017}
}


@article{badmus1996nonlinear,
  title={Nonlinear control of surge in axial compression systems},
  author={Badmus, OO and Chowdhury, S and Nett, CN},
  journal={Automatica},
  volume={32},
  number={1},
  pages={59--70},
  year={1996},
  publisher={Elsevier}
}

@article{gourdain2014large,
  title={Large eddy simulation of flows in industrial compressors: a path from 2015 to 2035},
  author={Gourdain, Nicolas and Sicot, Fr{\'e}d{\'e}ric and Duchaine, Florent and Gicquel, L},
  journal={Phil. Trans. of the Royal society A},
  volume={372},
  number={2022},
  pages={20130323},
  year={2014},
  publisher={The Royal Society Publishing}
}

@article{xiao2000center,
  title={Center Manifold of the Viscous Moore--Greitzer PDE Model},
  author={Xiao, Mingqing and Basar, Tamer},
  journal={SIAM Journal on Applied Mathematics},
  volume={61},
  number={3},
  pages={855--869},
  year={2000},
  publisher={SIAM}
}

@article{xiao2008quantitative,
  title={Quantitative characteristic of rotating stall and surge for Moore--Greitzer PDE model of an axial flow compressor},
  author={Xiao, MingQing},
  journal={SIAM Journal on Applied Dynamical Systems},
  volume={7},
  number={1},
  pages={39--62},
  year={2008},
  publisher={SIAM}
}



@article{liu2020smooth,
  title={Smooth Converse Lyapunov-Barrier Theorems for Asymptotic Stability with Safety Constraints and Reach-Avoid-Stay Specifications},
  author={Liu, Jun and Meng, Yiming and Li, Yinan and Fitzsimmons, Maxwell},
  journal={arXiv preprint arXiv:2009.04432},
  year={2020}
}

@book{belta2017formal,
  title={Formal Methods for Discrete-time Dynamical Systems},
  author={Belta, Calin and Yordanov, Boyan and Gol, Ebru Aydin},
  volume={89},
  year={2017},
  publisher={Springer}
}

@article{Girard2016multiscale,
author = {Girard, Antoine and Gossler, Gregor and Mouelhi, Sebti},
doi = {10.1109/TAC.2015.2478131},
journal = {IEEE Transactions on Automatic Control},
month = {jun},
number = {6},
pages = {1537--1549},
title = {Safety Controller Synthesis for Incrementally Stable Switched Systems Using Multiscale Symbolic Models},
volume = {61},
year = {2016}
}

@inproceedings{Hsu2018multilayer,
author = {Hsu, Kyle and Majumdar, Rupak and Mallik, Kaushik and Schmuck, Anne-Kathrin},
booktitle = {Proc. of HSCC},
pages = {120--129},
title = {Multi-Layered Abstraction-Based Controller Synthesis for Continuous-Time Systems},
year = {2018}
}

@book{filippov1988differential,
  title={Differential Equations with Discontinuous Righthand Sides},
  author={Filippov, Aleksei Fedorovich},
  year={1988},
  publisher={Springer Science \& Business Media}
}


@inproceedings{scots,
	Author = {Rungger, Matthias and Zamani, Majid},
	Booktitle = {Proc. of HSCC},
	Pages = {99-104},
	Title = {{SCOTS}: a Tool for the Synthesis of Symbolic Controllers},
	Year = {2016},
	}

@inproceedings{rocs,
	Author = {Li, Yinan and Liu, Jun},
	Booktitle = {Proc. of HSCC},
	Pages = {130--135},
	Title = {{ROCS}: A Robustly Complete Control Synthesis Tool for Nonlinear Dynamical Systems},
	Year = {2018}
	}


@article{kisielewicz2007stochastic,
  title={Stochastic differential inclusions and diffusion processes},
  author={Kisielewicz, Micha{\l}},
  journal={Journal of mathematical analysis and applications},
  volume={334},
  number={2},
  pages={1039--1054},
  year={2007},
  publisher={Elsevier}
}

@book{kisielewicz2013stochastic,
  title={Stochastic differential inclusions and applications},
  author={Kisielewicz, Micha{\l} and others},
  year={2013},
  publisher={Springer}
}

@article{lahijanian2015formal,
  title={Formal verification and synthesis for discrete-time stochastic systems},
  author={Lahijanian, Morteza and Andersson, Sean B and Belta, Calin},
  journal={IEEE Transactions on Automatic Control},
  volume={60},
  number={8},
  pages={2031--2045},
  year={2015},
  publisher={IEEE}
}

@inproceedings{cauchi2019efficiency,
  title={Efficiency through uncertainty: Scalable formal synthesis for stochastic hybrid systems},
  author={Cauchi, Nathalie and Laurenti, Luca and Lahijanian, Morteza and Abate, Alessandro and Kwiatkowska, Marta and Cardelli, Luca},
  booktitle={Proceedings of the 22nd ACM International Conference on Hybrid Systems: Computation and Control},
  pages={240--251},
  year={2019}
}

@inproceedings{majumdar2020symbolic,
  title={Symbolic controller synthesis for B{\"u}chi specifications on stochastic systems},
  author={Majumdar, Rupak and Mallik, Kaushik and Soudjani, Sadegh},
  booktitle={Proceedings of the 23rd International Conference on Hybrid Systems: Computation and Control},
  pages={1--11},
  year={2020}
}

@phdthesis{dutreix2020verification,
  title={Verification and synthesis for stochastic systems with temporal logic specifications},
  author={Dutreix, Maxence Dominique Henri},
  year={2020},
  school={Georgia Institute of Technology}
}


@article{dutreix2020specification,
  title={Specification-guided verification and abstraction refinement of mixed monotone stochastic systems},
  author={Dutreix, Maxence and Coogan, Samuel},
  journal={IEEE Transactions on Automatic Control},
  year={2020},
  publisher={IEEE}
}

@article{reissig2018symbolic,
  title={Symbolic optimal control},
  author={Reissig, Gunther and Rungger, Matthias},
  journal={IEEE Transactions on Automatic Control},
  volume={64},
  number={6},
  pages={2224--2239},
  year={2018},
  publisher={IEEE}
}

@article{vinod2017scalable,
  title={Scalable underapproximation for the stochastic reach-avoid problem for high-dimensional LTI systems using Fourier transforms},
  author={Vinod, Abraham P and Oishi, Meeko MK},
  journal={IEEE control systems letters},
  volume={1},
  number={2},
  pages={316--321},
  year={2017},
  publisher={IEEE}
}

@inproceedings{kariotoglou2013approximate,
  title={Approximate dynamic programming for stochastic reachability},
  author={Kariotoglou, Nikolaos and Summers, Sean and Summers, Tyler and Kamgarpour, Maryam and Lygeros, John},
  booktitle={2013 European Control Conference (ECC)},
  pages={584--589},
  year={2013},
  organization={IEEE}
}

@article{esfahani2016stochastic,
  title={The stochastic reach-avoid problem and set characterization for diffusions},
  author={Esfahani, Peyman Mohajerin and Chatterjee, Debasish and Lygeros, John},
  journal={Automatica},
  volume={70},
  pages={43--56},
  year={2016},
  publisher={Elsevier}
}

@book{subbaraman2015robust,
  title={Robust stability theory for stochastic dynamical systems},
  author={Subbaraman, Anantharaman},
  year={2015},
  publisher={University of California, Santa Barbara}
}

@inproceedings{poveda2015flexible,
  title={Flexible Nash seeking using stochastic difference inclusions},
  author={Poveda, Jorge I and Teel, Andrew R and Ne{\v{s}}i{\'c}, Dragan},
  booktitle={2015 American Control Conference (ACC)},
  pages={2236--2241},
  year={2015},
  organization={IEEE}
}

@article{malinowski2013interrelation,
  title={The interrelation between stochastic differential inclusions and set-valued stochastic differential equations},
  author={Malinowski, Marek T and Michta, Mariusz},
  journal={Journal of Mathematical Analysis and Applications},
  volume={408},
  number={2},
  pages={733--743},
  year={2013},
  publisher={Elsevier}
}

@book{mao2007stochastic,
  title={Stochastic differential equations and applications},
  author={Mao, Xuerong},
  year={2007},
  publisher={Elsevier}
}

@techreport{kushner1967stochastic,
  title={Stochastic stability and control},
  author={Kushner, Harold J},
  year={1967},
  institution={Brown Univ Providence RI}
}



@INPROCEEDINGS{meng2021control,

  author={Meng, Yiming and Li, Yinan and Liu, Jun},

  booktitle={2021 American Control Conference (ACC)}, 

  title={Control of Nonlinear Systems with Reach-Avoid-Stay Specifications: A Lyapunov-Barrier Approach with an Application to the Moore-Greizer Model}, 

  year={2021},

  volume={},

  number={},

  pages={2284-2291},

  doi={10.23919/ACC50511.2021.9483376}}



@inproceedings{wang2021barrier,
  title={Safety-Critical  Control  of  Stochastic  Systems  usingStochastic  Control  Barrier  Functions},
  author={Chuanzheng Wang$^\dagger$ and Yiming Meng$^\dagger$  and Stephen L. Smith and Jun Liu},
  booktitle={Proceedings of the Conference on Decisions and Control},
  year={2021}
}

@inproceedings{clark2019control,
  title={Control barrier functions for complete and incomplete information stochastic systems},
  author={Clark, Andrew},
  booktitle={2019 American Control Conference (ACC)},
  pages={2928--2935},
  year={2019},
  organization={IEEE}
}

@incollection{arnold1995random,
  title={Random dynamical systems},
  author={Arnold, Ludwig},
  booktitle={Dynamical systems},
  pages={1--43},
  year={1995},
  publisher={Springer}
}

@book{ethier2009markov,
  title={Markov processes: characterization and convergence},
  author={Ethier, Stewart N and Kurtz, Thomas G},
  volume={282},
  year={2009},
  publisher={John Wiley \& Sons}
}

@book{da2014stochastic,
  title={Stochastic equations in infinite dimensions},
  author={Da Prato, Giuseppe and Zabczyk, Jerzy},
  year={2014},
  publisher={Cambridge university press}
}

@article{prajna2007framework,
  title={A framework for worst-case and stochastic safety verification using barrier certificates},
  author={Prajna, Stephen and Jadbabaie, Ali and Pappas, George J},
  journal={IEEE Transactions on Automatic Control},
  volume={52},
  number={8},
  pages={1415--1428},
  year={2007},
  publisher={IEEE}
}

\begin{thebibliography}{10}

\bibitem{ames2016control}
Aaron~D Ames, Xiangru Xu, Jessy~W Grizzle, and Paulo Tabuada.
\newblock Control barrier function based quadratic programs for safety critical
  systems.
\newblock {\em IEEE Transactions on Automatic Control}, 62(8):3861--3876, 2016.

\bibitem{arnold1995random}
Ludwig Arnold.
\newblock Random dynamical systems.
\newblock In {\em Dynamical systems}, pages 1--43. Springer, 1995.

\bibitem{belta2017formal}
Calin Belta, Boyan Yordanov, and Ebru~Aydin Gol.
\newblock {\em Formal Methods for Discrete-time Dynamical Systems}, volume~89.
\newblock Springer, 2017.

\bibitem{cauchi2019efficiency}
Nathalie Cauchi, Luca Laurenti, Morteza Lahijanian, Alessandro Abate, Marta
  Kwiatkowska, and Luca Cardelli.
\newblock Efficiency through uncertainty: Scalable formal synthesis for
  stochastic hybrid systems.
\newblock In {\em Proceedings of the 22nd ACM International Conference on
  Hybrid Systems: Computation and Control}, pages 240--251, 2019.

\bibitem{clark2019control}
Andrew Clark.
\newblock Control barrier functions for complete and incomplete information
  stochastic systems.
\newblock In {\em 2019 American Control Conference (ACC)}, pages 2928--2935.
  IEEE, 2019.

\bibitem{da2014stochastic}
Giuseppe Da~Prato and Jerzy Zabczyk.
\newblock {\em Stochastic equations in infinite dimensions}.
\newblock Cambridge university press, 2014.

\bibitem{dutreix2020specification}
Maxence Dutreix and Samuel Coogan.
\newblock Specification-guided verification and abstraction refinement of mixed
  monotone stochastic systems.
\newblock {\em IEEE Transactions on Automatic Control}, 2020.

\bibitem{dutreix2020verification}
Maxence Dominique~Henri Dutreix.
\newblock {\em Verification and synthesis for stochastic systems with temporal
  logic specifications}.
\newblock PhD thesis, Georgia Institute of Technology, 2020.

\bibitem{esfahani2016stochastic}
Peyman~Mohajerin Esfahani, Debasish Chatterjee, and John Lygeros.
\newblock The stochastic reach-avoid problem and set characterization for
  diffusions.
\newblock {\em Automatica}, 70:43--56, 2016.

\bibitem{ethier2009markov}
Stewart~N Ethier and Thomas~G Kurtz.
\newblock {\em Markov processes: characterization and convergence}, volume 282.
\newblock John Wiley \& Sons, 2009.

\bibitem{faulwasser2009model}
Timm Faulwasser, Benjamin Kern, and Rolf Findeisen.
\newblock Model predictive path-following for constrained nonlinear systems.
\newblock In {\em Proc. of CDC}, pages 8642--8647. IEEE, 2009.

\bibitem{fribourg2013control}
Laurent Fribourg and Romain Soulat.
\newblock {\em Control of switching systems by invariance analysis: applcation
  to power electronics}.
\newblock John Wiley \& Sons, 2013.

\bibitem{girard2009approximately}
Antoine Girard, Giordano Pola, and Paulo Tabuada.
\newblock Approximately bisimilar symbolic models for incrementally stable
  switched systems.
\newblock {\em IEEE Transactions on Automatic Control}, 55(1):116--126, 2009.

\bibitem{kariotoglou2013approximate}
Nikolaos Kariotoglou, Sean Summers, Tyler Summers, Maryam Kamgarpour, and John
  Lygeros.
\newblock Approximate dynamic programming for stochastic reachability.
\newblock In {\em 2013 European Control Conference (ECC)}, pages 584--589.
  IEEE, 2013.

\bibitem{kisielewicz2007stochastic}
Micha{\l} Kisielewicz.
\newblock Stochastic differential inclusions and diffusion processes.
\newblock {\em Journal of mathematical analysis and applications},
  334(2):1039--1054, 2007.

\bibitem{kisielewicz2013stochastic}
Micha{\l} Kisielewicz et~al.
\newblock {\em Stochastic differential inclusions and applications}.
\newblock Springer, 2013.

\bibitem{kushner1967stochastic}
Harold~J Kushner.
\newblock Stochastic stability and control.
\newblock Technical report, Brown Univ Providence RI, 1967.

\bibitem{lahijanian2015formal}
Morteza Lahijanian, Sean~B Andersson, and Calin Belta.
\newblock Formal verification and synthesis for discrete-time stochastic
  systems.
\newblock {\em IEEE Transactions on Automatic Control}, 60(8):2031--2045, 2015.

\bibitem{li2020robustly}
Yinan Li and Jun Liu.
\newblock Robustly complete synthesis of memoryless controllers for nonlinear
  systems with reach-and-stay specifications.
\newblock {\em IEEE Transactions on Automatic Control}, 2020.

\bibitem{lin1996smooth}
Yuandan Lin, Eduardo~D Sontag, and Yuan Wang.
\newblock A smooth converse lyapunov theorem for robust stability.
\newblock {\em SIAM Journal on Control and Optimization}, 34(1):124--160, 1996.

\bibitem{liu2017robust}
Jun Liu.
\newblock Robust abstractions for control synthesis: Completeness via
  robustness for linear-time properties.
\newblock In {\em Proc. of HSCC}, pages 101--110, 2017.

\bibitem{liu2020converse}
Jun Liu.
\newblock Converse barrier functions via {L}yapunov functions.
\newblock {\em IEEE Transactions on Automatic Control}, 2021.

\bibitem{liu2020smooth}
Jun Liu, Yiming Meng, Yinan Li, and Maxwell Fitzsimmons.
\newblock Smooth converse lyapunov-barrier theorems for asymptotic stability
  with safety constraints and reach-avoid-stay specifications.
\newblock {\em arXiv preprint arXiv:2009.04432}, 2020.

\bibitem{majumdar2020symbolic}
Rupak Majumdar, Kaushik Mallik, and Sadegh Soudjani.
\newblock Symbolic controller synthesis for b{\"u}chi specifications on
  stochastic systems.
\newblock In {\em Proceedings of the 23rd International Conference on Hybrid
  Systems: Computation and Control}, pages 1--11, 2020.

\bibitem{malinowski2013interrelation}
Marek~T Malinowski and Mariusz Michta.
\newblock The interrelation between stochastic differential inclusions and
  set-valued stochastic differential equations.
\newblock {\em Journal of Mathematical Analysis and Applications},
  408(2):733--743, 2013.

\bibitem{mao2007stochastic}
Xuerong Mao.
\newblock {\em Stochastic differential equations and applications}.
\newblock Elsevier, 2007.

\bibitem{meng2021control}
Yiming Meng, Yinan Li, and Jun Liu.
\newblock Control of nonlinear systems with reach-avoid-stay specifications: A
  lyapunov-barrier approach with an application to the moore-greizer model.
\newblock In {\em 2021 American Control Conference (ACC)}, pages 2284--2291,
  2021.

\bibitem{nilsson2017augmented}
Petter Nilsson, Necmiye Ozay, and Jun Liu.
\newblock Augmented finite transition systems as abstractions for control
  synthesis.
\newblock {\em Discrete Event Dynamic Systems}, 27(2):301--340, 2017.

\bibitem{poveda2015flexible}
Jorge~I Poveda, Andrew~R Teel, and Dragan Ne{\v{s}}i{\'c}.
\newblock Flexible nash seeking using stochastic difference inclusions.
\newblock In {\em 2015 American Control Conference (ACC)}, pages 2236--2241.
  IEEE, 2015.

\bibitem{prajna2007framework}
Stephen Prajna, Ali Jadbabaie, and George~J Pappas.
\newblock A framework for worst-case and stochastic safety verification using
  barrier certificates.
\newblock {\em IEEE Transactions on Automatic Control}, 52(8):1415--1428, 2007.

\bibitem{reissig2018symbolic}
Gunther Reissig and Matthias Rungger.
\newblock Symbolic optimal control.
\newblock {\em IEEE Transactions on Automatic Control}, 64(6):2224--2239, 2018.

\bibitem{romdlony2016stabilization}
Muhammad~Zakiyullah Romdlony and Bayu Jayawardhana.
\newblock Stabilization with guaranteed safety using control lyapunov--barrier
  function.
\newblock {\em Automatica}, 66:39--47, 2016.

\bibitem{subbaraman2015robust}
Anantharaman Subbaraman.
\newblock {\em Robust stability theory for stochastic dynamical systems}.
\newblock University of California, Santa Barbara, 2015.

\bibitem{teel2000smooth}
Andrew~R Teel and Laurent Praly.
\newblock A smooth lyapunov function from a class-$\mathcal{KL}$ estimate
  involving two positive semidefinite functions.
\newblock {\em ESAIM: Control, Optimisation and Calculus of Variations},
  5:313--367, 2000.

\bibitem{vinod2017scalable}
Abraham~P Vinod and Meeko~MK Oishi.
\newblock Scalable underapproximation for the stochastic reach-avoid problem
  for high-dimensional lti systems using fourier transforms.
\newblock {\em IEEE control systems letters}, 1(2):316--321, 2017.

\bibitem{wang2021barrier}
Chuanzheng Wang$^\dagger$, Yiming Meng$^\dagger$, Stephen~L. Smith, and Jun
  Liu.
\newblock Safety-critical control of stochastic systems usingstochastic control
  barrier functions.
\newblock In {\em Proceedings of the Conference on Decisions and Control},
  2021.

\end{thebibliography}


\appendix

\noindent\textbf{Proof of Lemma \ref{lem: exist}. } 
The proof for non-emptiness and probability-$p$ invariance of $A$ is similar to \cite[Lemma 15]{liu2020smooth}, we can show that the reachable set within (ramdom) time interval $[\gamma_1, \gamma_2]$ is a valid choice by a  strong Markov property argument, where  $\gamma_1:=\inf\{t\geq 0: X_t\in\Gamma\}$ and $\gamma_2:=\inf\{t>\gamma_1: X_t\in\Gamma^c\}$.

Indeed, one can easily show that the reachable set $\bigcup_{X\in\bigcup_{x\in\init}\soln(x,W)}\mathscr{R}_\delta^{\gamma_1\leq t\leq \gamma_2}(X)\subset\Gamma$, and by the strong Markov property, for every $X\in\bigcup_{x\in\init}\soln(x,W)$, any restarted solution, denoted by $\tilde{X}\in\soln(X_{\gamma_3},W)$) where $\gamma_3\in[\gamma_1,\gamma_2)$,  has the same law as $X_{\gamma_3+s}$ for all $s\geq 0$. Continuing the above, for all $X\in \bigcup_{x\in\init}\soln(x,W)$, 
\begin{equation}\label{E: temp}
    \begin{split}
     &\ppp[\tilde{X}_{t\cj\sigma}\in\Gamma,\;\forall t\geq 0]\\
        \geq  & \inf\limits_{X\in\bigcup_{x\in\init}\soln(x,W)}\ppp[\gamma_1<\infty\;\text{and}\;\gamma_2=\infty]\geq p.
    \end{split}
\end{equation}
 By \eqref{E: temp}, 
  $\bigcup_{X\in\soln(\init,W)}\mathscr{R}_\delta^{\gamma_1\leq t\leq \gamma_2}(X)\subset A$. The probability-$p$ invariance is again by a standard strong Markov property argument and the definition of $A$.

To show that $A$ is closed,  let $x_n$ be a sequence in $A$ such that $x_n\rightarrow x\in\Gamma$, one can suppose the opposite and there there exists some $X\in\soln(x,W)$ such that 
$$\ppp[X_{t\cj\sigma}\in\Gamma, \;\forall t\geq 0]<p\Leftrightarrow\ppp[\tau<\infty]\geq 1-p,  $$
where $\tau:=\inf\{t\geq 0: X_t\in\Gamma^c\}$. 
Due to the weak compactness of the solution, by Skorohod, there exists a probability space $(\tilde{\Omega}^\dagger, \tilde{\mathscr{F}}^\dagger, \{\tilde{\mathscr{F}}_t^\dagger\},\tilde{\pp}^\dagger)$ and a sequence of processes $\tilde{X}^\sigma$, $\{(\tilde{X}^n)^\sigma\}$ such that there laws are $\ppp$ and $\{\ppp^n\}$ respectively and 
$
    \lim\limits_{n\rightarrow\infty}(\tilde{X}^n)^\sigma=\tilde{X}^\sigma,\;\; \tilde{\pp}^\dagger-\text{a.s.}
$ on every $[0,T]$.  On each $\{\tau\leq T\}$, $\tilde{X}^n_\tau \rightarrow \tilde{X}_\tau\notin \Gamma$, and since $\Gamma^c$ is open, for $n$ sufficiently large, we have $\tilde{X}^n_\tau\notin\Gamma$. The above shows that $\tau\leq T\implies \;\exists t\in[0,T]\;\text{s.t.}X^n\notin \Gamma$ for all $T$. Therefore, for sufficiently large $n$, sending $T$ to infinity, we have $\ppp^n[\exists t\geq 0\;\text{s.t.}\; X^n_t\notin \Gamma]\geq 1-p, $
which violates the probabilistic invariance of $A$. Hence, $x\in A$. The boundedness of $A$ is from the compactness of $\Gamma$. 
\hfill\pfbox

\noindent\textbf{Proof of Lemma \ref{lem: construct}.}
The proof is similar to \cite[Lemma 15]{liu2020converse} except for the context of weak solutions. We construct 
$Z_{s\cj\sigma}=X_{s\cj\sigma}+\frac{s}{T}[Z_{T\cj\sigma}-X_{T\cj\sigma}+(X_0-Y_0)]+Y_0-X_0$ for all $s\in[0,T]$. 
For any text function $\phi\in C^\infty(\R^n)$, we define
processes
\begin{small}
\begin{equation}
\begin{split}
     M^\phi(t) &=\phi(Z_{t\cj\sigma})-\phi(Z_0)\\
     & -\int_0^{t\cj\sigma} \nabla\phi(Z_s)\cdot f(Z_s)+\frac{1}{2}\operatorname{Tr}\left[(gg^T)(Z_s)\cdot \phi_{xx}(Z_s)\right] ds
\end{split}
\end{equation}
\end{small}
\begin{small}
\begin{equation}
   \begin{split}
     N^\phi(t) &=\phi(X_{t\cj\sigma})-\phi(X_0)\\
     & -\int_0^{t\cj\sigma} \nabla\phi(X_s)\cdot f(X_s)+\frac{1}{2}\operatorname{Tr}\left[(gg^T)(X_s)\cdot \phi_{xx}(X_s)\right] ds
\end{split}
\end{equation}
\end{small}
as well as a  martingale
\begin{equation}
    \hat{M}^\phi(t)=\phi(X_{t\cj\sigma})-\phi(X_0)
     -\int_0^{t\cj\sigma} \LL_d\phi(X_s) ds
\end{equation}
One can show that $|M^\phi(t)-\hat{M}^\phi(t)|$ has a bound $B(T,r, \phi, \delta')$ based on the properties of $f,g,\phi$ and $Z^\sigma$, i.e., 

\begin{equation}
    \begin{split}
 |M^\phi(t)-\hat{M}^\phi(t)| & \leq |M^\phi(t)-N^\phi(t)|+|N^\phi(t)-\hat{M}^\phi(t)|\\
 & \leq (2C_1r+TC_2r)+T\sup_{x\in D}\nabla\phi(x)\cdot \delta',
    \end{split}
\end{equation}
where $C_1$ is generated due to the choice of $\phi$, $C_2$ is due to the properties of $f,g$ and $\phi$. To make $M^\phi$ a family of martingales under the family of laws of the stopped process of $\sys$, one also needs to guarantee that 
$|M^\phi(t)-\hat{M}^\phi(t)|\leq T\sup_{x\in D}\nabla\phi(x)\cdot \delta$ for all $t\in[0,T]$,
since any two martingales 
of the martingale problem (under the laws of the corresponding stopped processes of $\sys$)  should not be differed larger than the above bound. 
Feasible ranges of $r$ can be obtained based on the requirement
$$r<\frac{\tau\sup_{x\in D}\nabla\phi(x)(\delta-\delta')}{2C_1+\tau C_2}.$$
The process $Z$ then satisfies the requirement given the feasible $r$.
\pfbox

\noindent\textbf{Proof of Proposition \ref{prop: attract}}
We just show the sketch. 
Without loss of generality, we consider $\oball_\eps(A)\subset \Gamma$. 
Suppose the claim is not true, 
then there exists some $\eps>0$ such that for all $n>0$ there exists $x_n\in \init$, $X^n\in\solnp(x_n,W)$ such that
$
    \ppp^n[|X_t^n|_A> \eps,\;\exists t\geq n]>1-p.
$

We now show this leads to a contradiction. By the assumption, $\sysp$ also satisfies the reach-avoid-stay specification $(\init,\Gamma, U, p)$.  By a similar argument of weak compactness as in the proof of Lemma \ref{lem: exist}, there exists a sufficiently large $N$ such that for all $n\geq N$,  defining the $\Gamma$-entering time $\gamma:=\inf\{t\geq 0: X_t^n\in\Gamma\}$,  we are able to show that  
\begin{equation}\label{E: oppo2}
    \ppp^n[\exists t\geq n: \;X_{t+\gamma}^n\notin \ball_\eps(A)]\geq 1-p
\end{equation}

Let $\tau:=N$ in Lemma \ref{lem: construct}. Note that, by the construction in Lemma \ref{lem: construct},  $\inf_{x\in A}|x|_\Gamma=0$ and  we are able to find a process $Z$ with $Z_0\in A$ and $Z_0\in \ball_r(X^n_{\gamma\cj\sigma})$ a.s.. By Lemma \ref{lem: construct}, there exists a process $Z\in\soln(Z_0,W)$ such that $Z_{t\cj\sigma}$ and $X^n_{(t+\gamma)\cj\sigma}$ share the same law for all $t\geq n$. However, $Z_{t\cj\sigma}\in A$ for all $t\geq 0$ by the definition of $A$. Therefore, $\ppp[\exists t\geq n: Z_{t\cj\sigma}\notin A]\leq \ppp[\exists t\geq 0: Z_{t\cj\sigma}\notin A]<1-p,$ which contradicts \eqref{E: oppo2}.
\hfill\pfbox

\noindent\textbf{Proof of Lemma \ref{lem: reccur}}. We just show the sketch.
The proof falls in a similar procedure as in the proof of \cite[Theorem 2.7]{mao2007stochastic}. We define the first hitting times $\tau_1, \tau_2$ of $\oball_{r_1}(A)$ and $\oball^c_{r_2}(A)$, where $0<r_1<r_2\leq R/2$. By It\^{o}'s formula, for each $X\in\soln(x,W)$ with  $x\in\oball_\eta(A)$ ($\eta$ is to be selected),  we have 
\begin{small}
\begin{equation*}
    \begin{split}
        0 &\leq V(x)+\eee^X\int_0^{\tau_1\cj\tau_2\cj t}\LL_d V(X(s))ds\\
        & \leq \alpha_2(\eta)-\alpha_3(r_1)\eee^X[\tau_1\cj\tau_2\cj t]
    \end{split}
\end{equation*}
\end{small}
On the other hand, 
\begin{equation*}
\begin{split}
      &\eee^X[\tau_1\cj\tau_2\cj t]\\ &\geq\int_\Omega \mathds{1}_{\{\tau_1\cj\tau_2\geq t\}}\cdot(\tau_1(\omega)\cj\tau_2(\omega)\cj t)\; d\ppp^X(\omega)\\
      &=t\ppp^X[\tau_1\cj\tau_2\geq t].
\end{split}
\end{equation*}
Combining the above, we have
$
    \ppp^X[\tau_1\cj\tau_2\geq t]\leq \alpha_2(\eta)/t\alpha_3(r_1)
$ for each $t$,
which holds for all $X\in\soln(x,W)$. By this relation, we construct $T:=T(\eps, \eta, r_1)=2\alpha_2(\eta)/\eps\alpha_3(r_1)$ 
and see
$$
    \infx\ppp^X[\tau_1\cj\tau_2< T]\geq 1-\frac{\eps}{2}.
$$
Now, let $\eta=\eta(\eps,r_2)$ be selected according to Remark \ref{rem: epsdelta} based on the Pr-U.S. property, such that $\ppp^X[\tau_2=\infty]\geq 1-\eps/2$ whenever $|x|_A\leq \eta$. Therefore, for $|x|_A\leq \eta$, we have for all $X\in\soln(x,W)$,
\begin{equation}
    \begin{split}
       1-\frac{\eps}{2}&\leq \ppp^X[\tau_1\cj\tau_2<
T] \leq \ppp^X[\tau_1<T]+\ppp^X[\tau_2<T]\\&
\leq \ppp^X[\tau_1<T]+\frac{\eps}{2},  
    \end{split}
\end{equation}
which complete the proof. \pfbox

\end{document}